\documentclass[12pt]{article}
\usepackage{latexsym}
\usepackage{amsmath}
\usepackage{amssymb}
\usepackage{amsfonts}
\usepackage{graphics,graphicx}

\numberwithin{equation}{section}

\newtheorem{theorem}{Theorem}[section]
\newtheorem{lemma}[theorem]{Lemma}
\newtheorem{pr}{Proposition}[section]

\newtheorem{remark}[theorem]{Remark}

\newcommand{\eproof}{{\mbox{\ }~\hfill
\mbox{\large $\Box$} \par \vskip 10pt}}
\newcommand{\pf}{\noindent{\bf Proof}}
\newcommand{\eps}{\varepsilon}

\newcommand{\R}{\mathbb R}

\title{ Carleman estimate for second order elliptic equations with Lipschitz leading coefficients and jumps at an interface}
\author{M.~ Di Cristo\thanks{Politecnico di Milano, Italy. Email: michele.dicristo@polimi.it}
\qquad E.~Francini\thanks{Universit\`a di Firenze, Italy. Email: elisa.francini@unifi.it}
\qquad C-L.~Lin\thanks{National Cheng Kung University, Taiwan. Email: cllin2@mail.ncku.edu.tw}\\
S.~Vessella\thanks{Universit\`a di Firenze, Italy. Email: sergio.vessella@unifi.it}
\qquad J-N.~Wang\thanks{National Taiwan University, Taiwan. Email: jnwang@math.ntu.edu.tw}}

\date{}

\begin{document}
\maketitle

\begin{abstract}
In this paper we prove a  local Carleman estimate for second order elliptic equations with a general anisotropic Lipschitz coefficients having a jump at an interface. Our approach does not rely on the techniques of microlocal analysis. We make use of the elementary method so that we are able to impose almost optimal assumptions on the coefficients and, consequently, the interface.
\end{abstract}

\tableofcontents

\section{Introduction}\label{intro}
\setcounter{equation}{0}

Since T. Carleman's pioneer work \cite{ca}, Carleman estimates have been indispensable tools for proving the unique continuation property for partial differential equations. Recently,  Carleman estimates have been successfully applied to study inverse problems, see for \cite{isakov}, \cite{KSU}.  Most of Carleman estimates are proved under the assumption that the leading coefficients possess certain regularity. For example, for general second order elliptic operators, Carleman estimates were proved when the leading coefficients are at least Lipschitz \cite{ho}, \cite{ho3}.  The restriction of regularity on the leading coefficients also reflects the fact that the unique continuation may fail if the coefficients are only H\"older continuous in $\R^n$ with $n\ge 3$ (see examples constructed by Pli\'s \cite{pl} and \cite{mi}). In $\R^2$, the unique continuation property holds for $W^{1,2}$ solutions of second elliptic equations in either non-divergence or divergence forms with essentially bounded coefficients \cite{bjs}, \cite{bn}, \cite{am}, \cite{s}. It should be noted that the unique continuation property for the second order elliptic equations in the plane with essentially bounded coefficients is deduced from the theory of quasiregular mappings. No Carleman estimates are derived in this situation.

From discussions above, Carleman estimates for second order elliptic operators with general discontinuous coefficients are not likely to hold. However,  when the discontinuities occur as jumps at an interface with homogeneous or non-homogeneous transmission conditions, one can still derive useful Carleman estimates. This is the main theme of the paper. There are some excellent works on this subject. We mention several closely related papers including Le Rousseau-Robbiano \cite{lr1}, \cite{lr2}, and Le Rousseau-Lerner \cite{ll}.  For the development of the problem and other related results, we refer the reader to the papers cited above and references therein. Our result is close to that of \cite{ll}, where the elliptic coefficient is a general anisotropic matrix-valued function. To put our paper in perspective, we would like to point out that the interface is assumed to be a $C^\infty$ hypersurface in \cite{ll} and the coefficients are $C^\infty$ away from the interface. Here we prove a local Carleman estimate near a flat interface from which it is easy to obtain under a standard  change of coordinates a Carleman estimate for operator with leading coefficients which have a jump discontinuity at a $C^{1,1}$ interface and are Lipschitz continuous apart from such an interface (see Theorem \ref{thm8.2} for a precise statement). The approach in \cite{ll} is close to Calder\'on's seminal work on the uniqueness of Cauchy problem \cite{cal} as an application of singular integral operators (or pseudo-differential operators). Therefore, the regularity assumptions of \cite{ll} are due to the use of calculus of pseudo-differential operators and the microlocal analysis techniques.

The aim here is to derive the Carleman estimate using more elementary methods. Our approach does not rely on the techniques of microlocal analysis, but rather on the straightforward Fourier transform. Thus we are able to relax the regularity assumptions on the coefficients and the interface. We first consider the simple case where the coefficients depend only on the normal variable.  Taking advantage of the simple structure of coefficients, we are able to derive a Carleman estimate by elementary computations with the help of the Fourier transform on the tangential variables. To handle the general coefficients, we rely on some type of partition of unity. In Section \ref{pre} after the Theorem \ref{thm8.2} we give a more detailed outline of our proof.

\section{Notations and statement of the main theorem}\label{pre}
Define $H_{\pm}=\chi_{\mathbb{R}^n_{\pm}}$ where $\mathbb{R}^n_{\pm}=\{(x,y)\in \mathbb{R}^{n-1}\times\mathbb{R}|y\gtrless0\}$ and $\chi_{\mathbb{R}^n_{\pm}}$ is the characteristic function of $\mathbb{R}^n_{\pm}$.
Let us stress that for a vector $(x,y)$ of $\mathbb R^n$, we mean
$x=(x_1,\dots,x_{n-1})\in\mathbb R^{n-1}$ and $y\in\mathbb R$.
In places we will use equivalently the symbols $D$, $\nabla$, $\partial$
to denote the gradient of a function and we will add the index $x$ or $y$
to denote gradient in $\mathbb R^{n-1}$ and the derivative with respect to $y$
respectively. 

Let $u_\pm\in C^\infty(\R^n)$. We define
\begin{equation*}\label{1.030}
u=H_+u_++H_-u_-=\sum_\pm H_{\pm}u_{\pm},
\end{equation*}
hereafter, we denote $\sum_\pm a_\pm=a_++a_-$,  and for $\mathbb{R}^{n-1}\times\mathbb{R}$
\begin{equation}\label{7.1}
\mathcal{L}(x,y,\partial)u:=\sum_{\pm}H_{\pm}{\rm div}_{x,y}(A_{\pm}(x,y)\nabla_{x,y}u_{\pm}),
\end{equation}
where
\begin{equation}\label{7.2}
A_{\pm}(x,y)=\{a^{\pm}_{ij}(x,y)\}^n_{i,j=1},\quad x\in \mathbb{R}^{n-1},y\in \mathbb{R}
\end{equation}
is a Lipschitz symmetric matrix-valued function satisfying, for given constants $\lambda_0\in (0,1]$, $M_0>0$,
\begin{equation}\label{7.3}
\lambda_0|z|^2\leq A_{\pm}(x,y)z\cdot z\leq \lambda^{-1}_0|z|^2,\, \forall (x,y)\in \mathbb{R}^n,\,\forall\, z\in \mathbb{R}^n
\end{equation}
and
\begin{equation}\label{7.4}
|A_{\pm}(x',y')-A_{\pm}(x,y)|\leq M_0(|x'-x|+|y'-y|).
\end{equation}
We write
\begin{equation}\label{7.5}
h_0(x):=u_+(x,0)-u_-(x,0),\ \forall\, x\in \mathbb{R}^{n-1},
\end{equation}
\begin{equation}\label{7.6}
h_1(x):=A_+(x,0)\nabla_{x,y}u_+(x,0)\cdot \nu-A_-(x,0)\nabla_{x,y}u_-(x,0)\cdot \nu,\ \forall\, x\in \mathbb{R}^{n-1},
\end{equation}
where $\nu=-e_n$.

Let us now introduce  the weight function. Let $\varphi$ be
\begin{equation}\label{2.1}
\varphi(y)=
\begin{cases}
\begin{array}{l}
\varphi_+(y):=\alpha_+y+\beta y^2/2,\quad y\geq 0,\\
\varphi_-(y):=\alpha_-y+\beta y^2/2,\quad y< 0,
\end{array}
\end{cases}
\end{equation}
where $\alpha_+$, $\alpha_-$ and $\beta$ are positive numbers which will be determined later. In what follows we denote by $\varphi_{+}$ and $\varphi_{-}$ the restriction of the weight function $\varphi$ to $[0,+\infty)$ and to $(-\infty,0)$ respectively. We use similar notation for any other weight functions. For any $\varepsilon>0$ let $$\psi_{\varepsilon}(x,y):=\varphi(y)-\frac{\varepsilon}{2}|x|^2,$$ and let

\begin{equation}\label{wei}
\phi_{\delta}(x,y):=\psi_{\delta}\left(\delta^{-1}x,\delta^{-1}y\right).
\end{equation}

For a function $h\in L^2(\mathbb{R}^{n})$, we define
\begin{equation*}
\hat{h}(\xi,y)=\int_{\mathbb{R}^{n-1}}h(x,y)e^{-ix\cdot\xi}\,dx,\quad \xi\in \mathbb{R}^{n-1}.
\end{equation*}
As usual we denote by $H^{1/2}(\mathbb{R}^{n-1})$ the space of the functions $f\in L^2(\mathbb{R}^{n-1})$ satisfying
$$\int_{\mathbb{R}^{n-1}}|\xi||\hat{f}(\xi)|^2d\xi<\infty,$$
with the norm
\begin{equation}\label{semR}
||f||^2_{H^{1/2}(\mathbb{R}^{n-1})}=\int_{\mathbb{R}^{n-1}}(1+|\xi|^2)^{1/2}|\hat{f}(\xi)|^2d\xi.
\end{equation}
Moreover we define
$$[f]_{1/2,\mathbb{R}^{n-1}}=\left[\int_{\mathbb{R}^{n-1}}\int_{\mathbb{R}^{n-1}}\frac{|f(x)-f(y)|^2}{|x-y|^n}dydx\right]^{1/2},$$
and recall that there is a positive constant $C$, depending only on $n$, such that
\begin{equation*}
C^{-1}\int_{\mathbb{R}^{n-1}}|\xi||\hat{f}(\xi)|^2d\xi\leq[f]^2_{1/2,\mathbb{R}^{n-1}}\leq C\int_{\mathbb{R}^{n-1}}|\xi||\hat{f}(\xi)|^2d\xi,
\end{equation*}
so that the norm \eqref{semR} is equivalent to the norm $||f||_{L^2(\mathbb{R}^{n-1})}+[f]_{1/2,\mathbb{R}^{n-1}}$. We use the letters $C, C_0, C_1, \cdots$ to denote constants. The value of the constants may change from line to line, but it is always greater than $1$.

We will denote by $B_r(x)$ the $(n-1)$-ball centered at $x\in \mathbb{R}^{n-1}$ with radius $r>0$. Whenever $x=0$ we denote $B_r=B_r(0)$.
\bigskip

\begin{theorem}\label{thm8.2}
Let $u$ and $A_{\pm}(x,y)$ satisfy \eqref{7.1}-\eqref{7.6}. There exist $\alpha_+,\alpha_-,\beta, \delta_0, r_0$ and $C$ depending on $\lambda_0, M_0$ such that if $\delta\le\delta_0$ and $\tau\geq C$, then
\begin{equation}\label{8.24}
\begin{aligned}
&\sum_{\pm}\sum_{k=0}^2\tau^{3-2k}\int_{\mathbb{R}^n_{\pm}}|D^k{u}_{\pm}|^2e^{2\tau\phi_{\delta,\pm}(x,y)}dxdy+\sum_{\pm}\sum_{k=0}^1\tau^{3-2k}\int_{\mathbb{R}^{n-1}}|D^k{u}_{\pm}(x,0)|^2e^{2\phi_\delta(x,0)}dx\\
&+\sum_{\pm}\tau^2[e^{\tau\phi_\delta(\cdot,0)}u_{\pm}(\cdot,0)]^2_{1/2,\mathbb{R}^{n-1}}+\sum_{\pm}[D(e^{\tau\phi_{\delta,\pm}}u_{\pm})(\cdot,0)]^2_{1/2,\mathbb{R}^{n-1}}\\
\leq &C\left(\sum_{\pm}\int_{\mathbb{R}^n_{\pm}}|\mathcal{L}( x, y,\partial)(u_{\pm})|^2\,e^{2\tau\phi_{\delta,\pm}(x,y)}dxdy+[e^{\tau\phi_\delta(\cdot,0)}h_1]^2_{1/2,\mathbb{R}^{n-1}}\right.\\
&\left.+[\nabla_x(e^{\tau\phi_\delta}h_0)(\cdot,0)]^2_{1/2,\mathbb{R}^{n-1}}+\tau^{3}\int_{\mathbb{R}^{n-1}}|h_0|^2e^{2\tau\phi_\delta(x,0)}dx+\tau\int_{\mathbb{R}^{n-1}}|h_1|^2e^{2\tau\phi_\delta(x,0)}dx\right).
\end{aligned}
\end{equation}
where $u=H_+u_++H_-u_-$,  $u_{\pm}\in C^\infty(\mathbb{R}^{n})$ and ${\rm supp}\, u\subset B_{\delta/2}\times[-\delta r_0,\delta r_0]$, and $\phi_\delta$ is given by \eqref{wei}.
\end{theorem}

\begin{remark}
Estimate \eqref{8.24} is a local Carleman estimate near $x_n=0$. As mentioned in the Introduction, by an easy change of coordinates, one can derive a local Carleman estimate near a $C^{1,1}$ interface from \eqref{8.24}.
\end{remark}
\begin{remark}\label{rem1}
Let us point out that the level sets $$\left\{(x,y)\in B_{\delta/2}\times (-\delta r_0, \delta r_0)):\phi_{\delta}(x,y)=t\right\}$$  have approximately the shape of paraboloid and, in a neighborhood of $(0,0)$, $\partial_y \phi_{\delta}>0$ so that the gradient of $\phi$ points inward the halfspace $\mathbb{R}_{+}^n$. These features are crucial to derive from the Carleman estimate \eqref{8.24} a H\"{o}lder type smallness propagation estimate across the interface $\left\{(x,0):x\in\mathbb{R}^{n-1}\right\}$ for weak solutions to the transmission problem
\begin{equation}
\begin{cases}
\begin{array}{l}
\mathcal{L}(x,y,\partial)u=\sum_{\pm}H_{\pm}b_{\pm}\cdot\nabla_{x,y}u_{\pm}+c_{\pm}u_{\pm},\\
u_{+}(x,0)-u_{+}(x,0)=0,\\
A_+(x,0)\nabla_{x,y}u_+(x,0)\cdot \nu-A_-(x,0)\nabla_{x,y}u_-(x,0)\cdot \nu=0,
\end{array}
\end{cases}
\end{equation}
where $b_{\pm}\in L^{\infty}(\mathbb{R}^n,\mathbb{R}^n)$ and $c_{\pm}\in L^{\infty}(\mathbb{R}^n) $. More precisely if the error of observation of $u$ is known in an open set of $\mathbb{R}_{+}^n$, we can find a H\"{o}lder control of $u$ in a bounded set of $\mathbb{R}_{-}^n$.
For more details about such type of estimate we refer to  \cite[Sect. 3.1]{lr1}.
\end{remark}

The proof of Theorem \ref{thm8.2} is divided into two steps as follows.

\medskip

\textbf{Step 1.} We first consider the particular case of the leading matrices \eqref{7.2} independent of $x$ and  we prove (Theorem \ref{pr22}), for the corresponding operator $\mathcal{L}(y,\partial)$, a Carleman estimate with the weight function $\phi(x,y)=\varphi(y)+ s\gamma\cdot x$, where $s$ is a suitable small number and $\gamma$ is an arbitrary unit vector of $\mathbb{R}^{n-1}$. The features of the leading matrices and of the weight function $\phi$ allow to factorize the Fourier transform of the conjugate of the operator $\mathcal{L}(y,\partial)u$ with respect to $\phi$. So that we can follow, roughly speaking, at an elementary level the strategy of \cite{ll} for the operator $\mathcal{L}(y,\partial)$. Nevertheless such an estimate has only a prepatory character to prove Theorem \ref{thm8.2}, because, due to the particular feature of the weight $\phi$ (i.e. linear with respect to $x$), the Carleman estimate obtained in Theorem \ref{pr22} cannot yield to any kind of significant smallness propagation estimate across the interface.

\medskip

\textbf{Step 2.} In the second we adapt the method described in \cite[Ch. 4.1]{Tr} to an operator with jump discontinuity. More precisely, we localize the operator \eqref{7.1} with respect to the $x$ variable and we linearize the weight function, again with respect the $x$ variable, and by the Carleman estimate obtained in the Step 1 we derive some local Carleman estimates. Subsequently we put together such local estimates by mean of the unity partition introduced in  \cite{Tr}.

\section{Step 1 - A Carleman estimate for leading coefficients depending on $y$ only}\label{sec1}
\setcounter{equation}{0}
In this section we consider the simple case of the leading matrices \eqref{7.2} independent of $x$. Moreover, the weight function that we consider is linear with respect to $x$ variable, so that, as explained above, the Carleman estimates we get here are only preliminary to the one that we will get in the general case.

Assume that
\begin{equation}\label{1.1}
A_{\pm}(y)=\{a^{\pm}_{ij}(y)\}^n_{i,j=1}
\end{equation}
are symmetric matrix-valued functions satisfying \eqref{7.3} and \eqref{7.4}, i.e.,
\begin{equation}\label{1.2}
\lambda_0|z|^2\leq A_{\pm}(y)z\cdot z\leq \lambda^{-1}_0|z|^2,\, \forall y\in \mathbb{R},\,\forall\, z\in \mathbb{R}^n
\end{equation}
\begin{equation}\label{1.3}
|A_{\pm}(y')-A_{\pm}(y'')|\leq M_0|y'-y''|, \quad\forall\, y',y''\in \mathbb{R}.
\end{equation}
From \eqref{1.2}, we have
\begin{equation}\label{Np5}
a_{nn}^{\pm}(y)\ge\lambda_0\quad\forall\, y\in\R.
\end{equation}

In the present case the the differential operator \eqref{7.1} became
\begin{equation}\label{1.6}
\mathcal{L}(y,\partial)u:=\sum_{\pm}H_{\pm}{\rm div}_{x,y}(A_{\pm}(y)\nabla_{x,y}u_{\pm}),
\end{equation}
where $u=\sum_\pm H_{\pm}u_{\pm}$, $u_\pm\in C^\infty(\R^n)$

We also set, for any $s\in[0,1]$ and $\gamma\in\mathbb{R}^{n-1}$ with $|\gamma|\le 1$
\begin{equation}\label{2.2}
\phi(x,y)=\varphi(y)+ s\gamma\cdot x=H_+\phi_++H_-\phi_-,
\end{equation}
where $\varphi$ is defined in \eqref{2.1}.

Our aim here is to prove the following Carleman estimate.
\begin{theorem}\label{pr22}
There exist $\tau_0$, $s_0$, $r_0$, $C$ and $\beta_0$ depending only on $\lambda_0$, $M_0$, such that for $\tau\geq\tau_0$, $0<s\le s_0<1$,
and for every  $w=\sum_{\pm}H_\pm w_\pm$ with ${\rm supp}\, w\subset B_1\times [-r_0,r_0]$, we have that
\begin{equation}\label{5.16}
\begin{aligned}
&\sum_{\pm}\sum_{k=0}^2\tau^{3-2k}\int_{\mathbb{R}^n_{\pm}}|D^k{w}_{\pm}|^2e^{2\tau\phi_\pm}dxdy\\
&+\sum_{\pm}\sum_{k=0}^1\tau^{3-2k}\int_{\mathbb{R}^{n-1}}|D^k{w}(x,0)|^2e^{2\tau\phi(x,0)}dx
+\sum_{\pm}\tau^2[(e^{\tau\phi}w_{\pm})(\cdot,0)]^2_{1/2,\mathbb{R}^{n-1}}\\
&+\sum_{\pm}[\partial_y(e^{\tau\phi_{\pm}}w_{\pm})(\cdot,0)]^2_{1/2,\mathbb{R}^{n-1}}
+\sum_{\pm}[\nabla_x(e^{\tau\phi}w_{\pm})(\cdot,0)]^2_{1/2,\mathbb{R}^{n-1}}\\
\leq &C\left(\int_{\mathbb{R}_\pm^n}|\mathcal{L}(y,\partial)w|^2e^{2\tau\phi_\pm}dxdy
+[\nabla_x\big(e^{\tau\phi(\cdot,0)}(w_+(\cdot,0)-w_-(\cdot,0))\big)]^2_{1/2,\mathbb{R}^{n-1}}\right.\\
&+[e^{\tau\phi(\cdot,0)}(A_+(0)\nabla_{x,y} w_+(x,0)\cdot \nu-A_-(0)\nabla_{x,y} w_-(x,0)\cdot \nu)]^2_{1/2,\mathbb{R}^{n-1}}\\
&+\tau\int_{\mathbb{R}^{n-1}}e^{2\tau\phi(x,0)}|A_+(0)\nabla_{x,y} w_+(x,0)\cdot \nu-A_-(0)\nabla_{x,y} w_-(x,0)\cdot \nu|^2dx\\
&\left.+\tau^3\int_{\mathbb{R}^{n-1}}e^{2\tau\phi(x,0)}|w_+(x,0)-w_-(x,0)|^2dx\right),
\end{aligned}
\end{equation}
with $\beta\ge\beta_0$ and $\alpha_\pm$ properly chosen.
\end{theorem}

\subsection{Fourier transform of the conjugate operator and its factorization}\label{Fourier}

To proceed further, we introduce some operators and find their properties. We use the notation $\partial_j=\partial_{x_j}$ for $1\le j\le n-1$.
Let us denote $B_{\pm}(y)=\{b^{\pm}_{jk}(y)\}^{n-1}_{j,k=1},$ the symmetric matrix satisfying, for $z=(z_1,\cdots,z_{n-1},z_n)=:(z',z_n)$,
\begin{equation}\label{1.12}
B_{\pm}(y)z'\cdot z'=A_{\pm}(y) z\cdot z\mid_{z_n=-\sum_{j=1}^{n-1}\frac{a_{nj}^{\pm}(y)z_j}{a^{\pm}_{nn}(y)}}.
\end{equation}
In view of \eqref{1.2} we have
\begin{equation}\label{1.15}
\lambda_1|z'|^2\leq B_{\pm}(y)z'\cdot z'\leq \lambda^{-1}_1|z'|^2,\, \quad\forall\, y\in \mathbb{R},\, \forall\, z'\in \mathbb{R}^{n-1},
\end{equation}
$\lambda_1\leq\lambda_0$  depends only on $\lambda_0$.

Notice that
\begin{equation}\label{1.8}
b_{jk}^{\pm}(y)=a_{jk}^{\pm}(y)-\frac{a_{nj}^{\pm}(y)a_{nk}^{\pm}(y)}{a_{nn}^{\pm}(y)},\quad j,k=1,\cdots,n-1.
\end{equation}
We define the operator
\begin{equation}\label{1.7}
T_{\pm}(y,\partial_x)u_{\pm}:=\sum_{j=1}^{n-1}\frac{a_{nj}^{\pm}(y)}{a_{nn}^{\pm}(y)}\partial_ju_{\pm}.
\end{equation}
It is easy to show, by direct calculations (\cite{ll}), that
\begin{equation}\label{1.9}
{\rm div}_{x,y}\big(A_{\pm}(y)\nabla_{x,y}u_{\pm}\big)=(\partial_y+T_{\pm})a_{nn}^{\pm}(y)(\partial_y+T_{\pm})u_{\pm}+{\rm div}_{x}\big(B_{\pm}(y)\nabla_{x}u_{\pm}\big).
\end{equation}
Let $w=\sum_{\pm}H_{\pm}w_{\pm}$, where $w_{\pm}\in C^\infty(\mathbb{R}^{n})$, and define
\begin{equation}\label{1.4}
\theta_0(x):=w_+(x,0)-w_-(x,0)\quad \mbox{ for } x\in \mathbb{R}^{n-1},
\end{equation}
\begin{equation}\label{1.5}
\theta_1(x):=A_+(0)\nabla_{x,y}w_+(x,0)\cdot\nu-A_-(0)\nabla_{x,y}w_-(x,0)\cdot\nu\quad  \mbox{ for }  x\in \mathbb{R}^{n-1},
\end{equation}
where $\nu=-e_n$. By straightforward calculations, we get
\begin{equation}\label{1.19}
a_{nn}^{+}(y)\big(\partial_y+T_{+}(y,\partial_x)\big)w_{+}(x,y)\mid_{y=0}-a_{nn}^{-}(y)\big(\partial_y+T_{-}(y,\partial_x)\big)w_{-}(x,y)\mid_{y=0}=-\theta_1(x).
\end{equation}

 In order to derive the Carleman estimate \eqref{5.16}, we investigate the conjugate operator of $\mathcal{L}(y,\partial)$ with $e^{\tau\phi}$ for $\phi$ given by \eqref{2.2}. Let $v=e^{\tau\phi}w$ and $\tilde{v}=e^{-\tau s\gamma\cdot x}v$, then we have
$$w=e^{-\tau\phi}v=\sum_{\pm}H_{\pm}e^{-\tau\phi_\pm}v_{\pm}=\sum_{\pm}H_{\pm}e^{-\tau\varphi_\pm}\tilde{v}_{\pm}$$
and therefore
\[
e^{\tau\phi}\mathcal{L}(y,\partial)(e^{-\tau\phi}v)=e^{\tau s\gamma\cdot x}e^{\tau\varphi}\mathcal{L}(y,\partial)(e^{-\tau\varphi}\tilde{v}).
\]
It follows from \eqref{1.9}  that
\[
\begin{aligned}
e^{\tau\varphi}\mathcal{L}(y,\partial)(e^{-\tau\varphi}\tilde{v})
=&\sum_{\pm}H_{\pm}\big[(\partial_y-\tau\varphi_{\pm}'+T_{\pm})a_{nn}^{\pm}(y)(\partial_y-\tau\varphi_{\pm}'+T_{\pm})\tilde{v}_{\pm}\big]\\
&\,+\sum_{\pm}H_{\pm}{\rm div}_{x}\big(B_{\pm}(y)\nabla_{x}\tilde{v}_{\pm}\big),
\end{aligned}
\]
which leads to
\begin{equation}\label{2.5}
\begin{aligned}
&e^{\tau\phi}\mathcal{L}(y,\partial)(e^{-\tau\phi}v)
=e^{\tau s\gamma\cdot x}e^{\tau\varphi}\mathcal{L}(y,\partial)(e^{-\tau\varphi}\tilde{v})\\
&\,\,=e^{\tau s\gamma\cdot x}\sum_{\pm}H_{\pm}\big[(\partial_y-\tau\varphi_{\pm}'+T_{\pm})a_{nn}^{\pm}(y)(\partial_y-\tau\varphi_{\pm}'+T_{\pm})(e^{-\tau s\gamma\cdot x}v_{\pm})\big]\\
&\,\,\,+e^{\tau s\gamma\cdot x}\sum_{\pm}H_{\pm}{\rm div}_{x}\big(B_{\pm}(y)\nabla_{x}(e^{-\tau s\gamma\cdot x}v_{\pm})\big).
\end{aligned}
\end{equation}
By the definition of $T_{\pm}(y,\partial_x)$, we get that
\begin{eqnarray*}
T_{\pm}(y,\partial_x)(e^{-\tau s\gamma\cdot x}v_{\pm})&=&e^{-\tau s\gamma\cdot x}\sum_{j=1}^{n-1}\frac{a_{nj}^{\pm}(y)}{a_{nn}^{\pm}(y)}(\partial_jv_{\pm}-\tau s\gamma_jv_{\pm})\\
&:=&e^{-\tau s\gamma\cdot x}\,T_{\pm}(y,\partial_x-\tau s\gamma)v_{\pm}.
\end{eqnarray*}
To continue the computation, we observe that
\begin{equation}\label{2.7}
\begin{aligned}
&e^{\tau s\gamma\cdot x}\big[(\partial_y-\tau\varphi_{\pm}'+T_{\pm}(y,\partial_x))a_{nn}^{\pm}(y)(\partial_y-\tau\varphi_{\pm}'+T_{\pm}(y,\partial_x))(e^{-\tau s\gamma\cdot x}v_{\pm})\big]\\
=&\big(\partial_y-\tau\varphi_{\pm}'+T_{\pm}(y,\partial_x-\tau s\gamma)\big)a_{nn}^{\pm}(y)\big(\partial_y-\tau\varphi_{\pm}'+T_{\pm}(y,\partial_x-\tau s\gamma)\big)v_{\pm}
\end{aligned}
\end{equation}

and
\begin{equation}\label{2.8}
\begin{aligned}
&e^{\tau s\gamma\cdot x}{\rm div}_{x}\big(B_{\pm}(y)\nabla_{x}(e^{-\tau s\gamma\cdot x}v_{\pm})\big)\\
&\,\,=\sum_{j,k=1}^{n-1}b_{jk}^{\pm}(y)\partial^2_{jk}v_{\pm}-2s\tau\sum_{j,k=1}^{n-1}b_{jk}^{\pm}(y)\partial_{j}v_{\pm}\gamma_k
+s^2\tau^2\sum_{j,k=1}^{n-1}b_{jk}^{\pm}(y)\gamma_j\gamma_kv_{\pm}.
\end{aligned}
\end{equation}
Combining \eqref{2.5}, \eqref{2.7} and \eqref{2.8} yields
\begin{equation}\label{2.9}
\begin{aligned}
&e^{\tau\phi}\mathcal{L}(y,\partial)(e^{-\tau\phi}v)\\
=&\sum_{\pm}H_{\pm}\big[\big(\partial_y-\tau\varphi_{\pm}'+T_{\pm}(y,\partial_x-\tau s\gamma)\big)a_{nn}^{\pm}(y)\big(\partial_y-\tau\varphi_{\pm}'+T_{\pm}(y,\partial_x-\tau s\gamma)\big)v_{\pm}\big]\\
&+\sum_{\pm}H_{\pm}\big[\sum_{j,k=1}^{n-1}b_{jk}^{\pm}(y)\partial^2_{jk}v_{\pm}-2s\tau\sum_{j,k=1}^{n-1}b_{jk}^{\pm}(y)\partial_{j}v_{\pm}\gamma_k+s^2\tau^2\sum_{j,k=1}^{n-1}b_{jk}^{\pm}(y)\gamma_j\gamma_kv_{\pm}\big].
\end{aligned}
\end{equation}

We now focus on the analysis of $e^{\tau\phi}\mathcal{L}(y,\partial)(e^{-\tau\phi}v)$. To simplify it, we introduce some notations:
\begin{equation}\label{2.10}
f(x,y)=e^{\tau\phi}\mathcal{L}(y,\partial)(e^{-\tau\phi}v),
\end{equation}
\begin{equation}\label{2.11}
B_{\pm}(\xi,\gamma,y)=\sum_{j,k=1}^{n-1}b_{jk}^{\pm}(y)\xi_j\gamma_k,\quad \xi\in \mathbb{R}^{n-1},
\end{equation}
\begin{equation}\label{2.14}
\zeta_{\pm}(\xi,y)=\frac{1}{a_{nn}^{\pm}(y)}\big[B_{\pm}(\xi,\xi,y)+2is\tau B_{\pm}(\xi,\gamma,y)-s^2\tau^2B_{\pm}(\gamma,\gamma,y)\big],
\end{equation}
and
\begin{equation}\label{2.12}
t_{\pm}(\xi,y)=\sum_{j=1}^{n-1}\frac{a_{nj}^{\pm}(y)}{a_{nn}^{\pm}(y)}\xi_j.
\end{equation}
By \eqref{2.9}, we have
\begin{equation}\label{2.15}
\hat{f}(\xi,y)
=\sum_{\pm}H_{\pm}P_{\pm}\hat{v}_{\pm},
\end{equation}
where
\begin{equation}\label{2.17}
\begin{aligned}
P_\pm\hat{v}_\pm:=&\big(\partial_y-\tau\varphi_{\pm}^\prime+it_{\pm}(\xi+i\tau s\gamma,y)\big)a_{nn}^{\pm}(y)\big(\partial_y-\tau\varphi_{\pm}^\prime+it_{\pm}(\xi+i\tau s\gamma,y)\big)\hat{v}_{\pm}\\
&-a_{nn}^{\pm}(y)\zeta_{\pm}(\xi,y)\hat{v}_{\pm}.
\end{aligned}
\end{equation}

Our aim is to estimate  $f(x,y)$ or, equivalently, its Fourier transform  $\hat{f}(\xi,y)$.
In order to do this, we want to factorize the operators $P_{\pm}$. For any $z=a+ib$ with $(a,b)\neq(0,0)$, we define the square root of $z$,
$$\sqrt{z}=\sqrt{\frac{a+\sqrt{a^2+b^2}}{2}}+i\frac{b}{\sqrt{2(a+\sqrt{a^2+b^2})}}.$$
We remind that the square root $\sqrt{z}$ is defined with a cut along the negative real axis and note that $\Re \sqrt{z}\geq 0$. Thus, it needs extra care to estimate its derivative. Now we define two operators
\begin{equation}\label{3.10}
E_{\pm}=\partial_y+it_{\pm}(\xi+i\tau s\gamma,y)-(\tau\varphi^\prime_{\pm}+\sqrt{\zeta_{\pm}}),
\end{equation}
\begin{equation}\label{3.11}
F_{\pm}=\partial_y+it_{\pm}(\xi+i\tau s\gamma,y)-(\tau\varphi^\prime_{\pm}-\sqrt{\zeta_{\pm}}).
\end{equation}
With all the definitions given above, we obtain that
\begin{equation}\label{3.12}
P_+\hat{v}_+=E_{+}a_{nn}^+(y)F_+\hat{v}_+-\hat{v}_+\partial_y\big(a_{nn}^+(y)\sqrt{\zeta_{+}}\big),
\end{equation}
\begin{equation}\label{3.13}
P_-\hat{v}_-=F_{-}a_{nn}^-(y)E_-\hat{v}_-+\hat{v}_-\partial_y\big(a_{nn}^-(y)\sqrt{\zeta_{-}}\big).
\end{equation}

Let us now introduce some other useful notations and estimates that will be intensively used in the sequel. After
taking the Fourier transform, the terms on the interface \eqref{1.4} and  \eqref{1.19}, become
\begin{equation}\label{2.18}
\eta_0(\xi):=\hat{v}_+(\xi,0)-\hat{v}_-(\xi,0)=\widehat{e^{\tau\phi(x,0)}\theta_0(x)}
\end{equation}
and
\begin{equation}\label{2.19}
\begin{aligned}
&\eta_1(\xi):=\widehat{-e^{\tau\phi(x,0)}\theta_1(x)}\\
&\,=a_{nn}^{+}(0)\big[\partial_y\hat{v}_{+}(\xi,0)-\tau\alpha_+\hat{v}_+(\xi,0)+it_{+}(\xi+i\tau s\gamma,0)\hat{v}_+(\xi,0)\big]\\
&\,\,\,-a_{nn}^{-}(0)\big[\partial_y\hat{v}_{-}(\xi,0)-\tau\alpha_-\hat{v}_-(\xi,0)+it_{-}(\xi+i\tau s\gamma,0)\hat{v}_-(\xi,0)\big].
\end{aligned}
\end{equation}
For simplicity, we denote
\begin{equation}\label{2.20}
V_{\pm}(\xi):=a_{nn}^{\pm}(0)\big[\partial_y\hat{v}_{\pm}(\xi,0)-\tau\alpha_{\pm}\hat{v}_{\pm}(\xi,0)+it_{\pm}(\xi+i\tau s\gamma,0)\hat{v}_{\pm}(\xi,0)\big],
\end{equation}
so that
\begin{equation}\label{2.20.1}
	V_{+}(\xi)-V_{-}(\xi)=\eta_1(\xi).
\end{equation}
Moreover, we define
\begin{equation*}
m_{\pm}(\xi,y):=\sqrt{\frac{B_{\pm}(\xi,\xi,y)}{a_{nn}^{\pm}(y)}}.
\end{equation*}
From \eqref{1.15} 
we have
\begin{equation}\label{3.1}
\lambda_1|\xi|^2\leq B_{\pm}(\xi,\xi,y)\leq \lambda^{-1}_1|\xi|^2,\quad\forall\, y\in\R,\, \forall\,\xi\in\R^{n-1},
\end{equation}
and, from \eqref{1.3},
\begin{equation}\label{3.2}
|\partial_yB_{\pm}(\xi,\eta,y)|\leq M_1|\xi||\eta|,\quad \forall\, \xi,\,\eta\in \mathbb{R}^{n-1},
\end{equation}
where $M_1$ depends only on  $\lambda_0$ and $M_0$.
In a similar way, we list here some useful bounds, that can be easily obtained from \eqref{1.15} and \eqref{1.3}.
\begin{equation}\label{3.6}
\lambda_2|\xi|\leq m_{\pm}(\xi,y)\leq \lambda^{-1}_2|\xi|,
\end{equation}
\begin{equation}\label{3.7}
|\partial_ym_{\pm}(\xi,y)|\leq M_2|\xi|,
\end{equation}
\begin{equation}\label{3.6.1}
|t_{\pm}(\xi,y)|\leq \lambda^{-1}_3|\xi|,
\end{equation}
\begin{equation}\label{3.7.1}
|\partial_yt_{\pm}(\xi,y)|\leq M_3|\xi|,
\end{equation}
\begin{equation}\label{3.8}
|\zeta_{\pm}(\xi,y)|\leq (\lambda_0\lambda_1)^{-1}(|\xi|^2+s^2\tau^2),
\end{equation}
\begin{equation}\label{3.9}
|\partial_y\zeta_{\pm}(\xi,y)|\leq M_4(|\xi|^2+s^2\tau^2).
\end{equation}
Here  $\lambda_2=\sqrt{\lambda_0\lambda_1}$, $\lambda_3$ depends only on $\lambda_0$, while $M_2$, $M_3$ and $M_4$ depends only on $\lambda_0$ and $M_0$.

\subsection{Derivation of the Carleman estimate for the simple case}\label{simple}

The derivation of the Carleman estimate \eqref{5.16} is a simple consequence of the auxiliary Proposition \ref{pr22}  stated below 
and proved in the following Section \ref{keyprop} via the inverse Fourier transform. We first set
\begin{equation*}
L:=\sup_{\xi\in \mathbb{R}^{n-1}\setminus\{0\}}\frac{m_+(\xi,0)}{m_-(\xi,0)}.
\end{equation*}
Note that, by \eqref{3.6}, $\lambda_2^{2}\leq L\leq \lambda_2^{-2}$. Now we introduce the fundamental assumption on the coefficients $\alpha_{\pm}$ in the weight function. As in \cite{ll}, we choose positive $\alpha_+$ and $\alpha_-$, such that
\begin{equation}\label{4.3}
L<\frac{\alpha_+}{\alpha_-}.
\end{equation}
This choice will only be conditioned by $\lambda_0$. These constants will be fixed. Denote the factor
\[
\Lambda=(|\xi|^2+\tau^2)^{1/2}.
\]
We now state our main tool.
\begin{pr}\label{bigthm}
There exist $\tau_0$, $s_0$, $\rho$, $\beta$ and $C$, depending only on $\lambda_0$ and $M_0$, such that for $\tau\geq\tau_0$, ${\rm supp}\, \hat{v}_{\pm}(\xi,\cdot)\subset [-\rho,\rho]$, $s\leq s_0<1$, we have
\begin{eqnarray}\label{4.4}
&&\frac{1}{\tau}\sum_{\pm}||\partial^2_y\hat{v}_{\pm}(\xi,\cdot)||^2_{L^2(\mathbb{R}_\pm)}
+\frac{\Lambda^2}{\tau}\sum_{\pm}||\partial_y\hat{v}_{\pm}(\xi,\cdot)||^2_{L^2(\mathbb{R}_\pm)}
\nonumber\\
&&\quad+\frac{\Lambda^4}{\tau}\sum_{\pm}||\hat{v}_{\pm}(\xi,\cdot)||^2_{L^2(\mathbb{R}_\pm)}+\Lambda\sum_{\pm}|V_{\pm}(\xi)|^2+\Lambda^3\sum_{\pm}|\hat{v}_{\pm}(\xi,0)|^2\nonumber\\
 &&\quad\leq C\left(\sum_{\pm}||P_{\pm}\hat{v}_{\pm}(\xi,\cdot)||^2_{L^2(\mathbb{R}_\pm)}+\Lambda|\eta_1(\xi)|^2+\Lambda^3|\eta_0(\xi)|^2\right).
\end{eqnarray}
Here $\R_{\pm}=\{y\in\R\,:\,y\gtrless 0\}$.
\end{pr}

\bigskip\noindent
{\bf Proof of Theorem~\ref{pr22}}. Substituting \eqref{2.15} and the definitions of $\eta_0$, $\eta_1$ (see \eqref{2.18}, \eqref{2.19}) into the right hand side of \eqref{4.4} implies \begin{equation}\label{5.13}
\begin{aligned}
&\frac{1}{\tau}\sum_{\pm}||\partial^2_y\hat{v}_{\pm}(\xi,\cdot)||^2_{L^2(\mathbb{R}_{\pm})}
+\frac{\Lambda^2}{\tau}\sum_{\pm}||\partial_y\hat{v}_{\pm}(\xi,\cdot)||^2_{L^2(\mathbb{R}_{\pm})}
+\frac{\Lambda^4}{\tau}\sum_{\pm}||\hat{v}_{\pm}(\xi,\cdot)||^2_{L^2(\mathbb{R}_{\pm})}\\
&+\Lambda\sum_{\pm}|V_{\pm}(\xi)|^2+\Lambda^3\sum_{\pm}|\hat{v}_{\pm}(\xi,0)|^2\\
\leq &C\left(\sum_{\pm}||f(\xi,\cdot)||^2_{L^2(\mathbb{R})}+\Lambda|\widehat{e^{\tau\phi(\cdot,0)}\theta_1(\cdot)}|^2
+\Lambda^3|\widehat{e^{\tau\phi(\cdot,0)}\theta_0(\cdot)}|^2\right).
\end{aligned}
\end{equation}
Recalling \eqref{2.20}, 
it is not hard to see that
\begin{equation}\label{5.14}
\Lambda\sum_{\pm}|\partial_y\hat{v}_{\pm}(\xi,0)|^2\leq C\left(\Lambda\sum_{\pm}|V_{\pm}(\xi)|^2+\Lambda^3\sum_{\pm}|\hat{v}_{\pm}(\xi,0)|^2\right).
\end{equation}
Since $\Lambda^4\geq |\xi|^2\tau^2+|\xi|^4+\tau^4$, $|\xi|^3+|\xi|^2\tau+|\xi|\tau^2+\tau^3\le C\Lambda^3$, and $\Lambda^3\le C'(|\xi|^3+\tau^3)$, by integrating in $\xi$,
we can deduce from  \eqref{5.13} and \eqref{5.14} that
\begin{equation}\label{5.15}
\begin{aligned}
&\sum_{\pm}\sum_{k=0}^2\tau^{3-2k}\int_{\mathbb{R}^n_{\pm}}|D^k{v}_{\pm}|^2
+\sum_{\pm}[\nabla_xv_{\pm}(\cdot,0)]^2_{1/2,\mathbb{R}^{n-1}}
+\sum_{\pm}[\partial_yv_{\pm}(\cdot,0)]^2_{1/2,\mathbb{R}^{n-1}}\\
&+\sum_{\pm}\tau^2[v_{\pm}(\cdot,0)]^2_{1/2,\mathbb{R}^{n-1}}+\sum_{\pm}\tau\int_{\mathbb{R}^{n-1}}|\nabla_xv_{\pm}(x,0)|^2dx\\
&+\sum_{\pm}\tau\int_{\mathbb{R}^{n-1}}|\partial_yv_{\pm}(x,0)|^2dx+\sum_{\pm}\tau^3\int_{\mathbb{R}^{n-1}}|v_{\pm}(x,0)|^2dx\\
\leq &C\left(||f||^2_{L^2(\mathbb{R}^n)}+[e^{\tau\phi(\cdot,0)}\theta_1(\cdot)]^2_{1/2,\mathbb{R}^{n-1}}
+[\nabla_x\big(e^{\tau\phi(\cdot,0)}\theta_0(\cdot)\big)]^2_{1/2,\mathbb{R}^{n-1}}\right.\\
&\left.+\tau\int_{\mathbb{R}^{n-1}}e^{2\tau\phi(x,0)}|\theta_1|^2dx+\tau^3\int_{\mathbb{R}^{n-1}}e^{2\tau\phi(x,0)}|\theta_0|^2dx\right).
\end{aligned}
\end{equation}
Replacing $v_{\pm}=e^{\tau\phi_{\pm}}w_{\pm}$ into \eqref{5.15} immediately leads to \eqref{5.16}.\eproof

\subsection{Proof of Proposition \ref{pr22}}\label{keyprop}

Let $\kappa$ be the positive number
\begin{equation}\label{kappa}
	\kappa=\frac{1}{2}\left(1-L\frac{\alpha_-}{\alpha_+}\right)
\end{equation}
depending only on $\lambda_0$ and $M_0$. The proof of Proposition \ref{bigthm} will be divided into three cases
\begin{equation*}
\left\{
\begin{aligned}
&\tau\geq \frac{\lambda_2^2|\xi|}{2s_0},\\
&\frac{m_+(\xi,0)}{(1-\kappa)\,\alpha_+}\leq\tau\leq\frac{\lambda_2^2|\xi|}{2s_0},\\
&\tau\leq\frac{m_+(\xi,0)}{(1-\kappa)\,\alpha_+}.
\end{aligned}
\right.
\end{equation*}
Recall that $\lambda_2=\sqrt{\lambda_0\lambda_1}$  (from \eqref{3.6}) depends only on $\lambda_0$. Of course, we first choose a small $s_0<1$, depending on $\lambda_0$ and $M_0$ only, such that
\[
\frac{m_+(\xi,0)}{(1-\kappa)\,\alpha_+}\leq\frac{\lambda_2^2|\xi|}{2s_0},\quad\forall\,\xi\in\R^n.
\]
A smaller value $s_0$ will be chosen later in the proof.

We need to introduce here some further notations.
First of all, let us denote by
\[P^0_{\pm},\; E^0_{\pm},\; \mbox{and}\; F^0_{\pm}\]
the operators defined by \eqref{2.17}, \eqref{3.10} and \eqref{3.11}, respectively, in the special case $s=0$. We also give special names to these functions that will be used in the proof:
\begin{equation}\label{omega}
	\omega_+(\xi,y)=a_{nn}^+(y)F_+\hat{v}_+(\xi,y),\quad \omega_-(\xi,y)=a_{nn}^-(y)E_-\hat{v}_-(\xi,y)
\end{equation}
and, for the special case $s=0$,
\begin{equation}\label{omega0}
\omega^0_+(\xi,y)=a_{nn}^+(y)F^0_+\hat{v}_+(\xi,y),\quad	\omega^0_-(\xi,y)=a_{nn}^-(y)E^0_-\hat{v}_-(\xi,y).
\end{equation}

\bigskip
\noindent{\bf Case 1:}
\begin{equation}\label{case1}
	\tau\geq \frac{\lambda_2^2|\xi|}{2s_0}
\end{equation}
%
Note that, in this case, we have $|\xi|\leq 2\lambda_2^{-2} s_0\tau$, which implies
\begin{equation}\label{4.5}
\tau\leq\Lambda\leq \sqrt{5}\lambda_2^{-2}\tau.
\end{equation}

\noindent
We will need several lemmas. In the first one, we estimate the difference $P_{\pm}\hat{v}_{\pm}-P^0_{\pm}\hat{v}_{\pm}$.
\begin{lemma}\label{lem4.1}
Let $\tau\geq1$ and assume \eqref{case1}, then we have
\begin{equation}\label{4.7}
|P_{\pm}\hat{v}_{\pm}(\xi,y)-P^0_{\pm}\hat{v}_{\pm}(\xi,y)|
\leq C s\tau\big[\tau(\alpha_{\pm}+1+\beta|y|)|\hat{v}_{\pm}(\xi,y)|+|\partial_y\hat{v}_{\pm}(\xi,y)|\big],
\end{equation}
where $C$ depends only on $\lambda_0$ and $M_0$.
\end{lemma}
\pf. It should be noted that
\[
\zeta_{\pm}(\xi,y)|_{s=0}=\frac{B_\pm(\xi,\xi,y)}{a_{nn}^{\pm}(y)}.
\]
By simple calculations, and dropping $\pm$ for the sake of shortness, we can write
\begin{equation}\label{4.7.1}
	P\hat{v}(\xi,y)-P^0\hat{v}(\xi,y)=I_1+I_2+I_3,
\end{equation}
where
\begin{eqnarray*}
I_1&=&\big(it(\xi+i\tau s\gamma,y)-it(\xi,y)\big)a_{nn}(y)\big(\partial_y-\tau\varphi'+it(\xi+i\tau s\gamma,y)\big)\hat{v},\\
I_2&=&\big(\partial_y-\tau\varphi'+it(\xi,y)\big)a_{nn}(y)\big(it(\xi+i\tau s\gamma,y)-it(\xi,y)\big)\hat{v},
\end{eqnarray*}
and
\begin{equation*}
	I_3=a_{nn}^{\pm}(y)\zeta_{\pm}(\xi,y)-B_{\pm}(\xi,\xi,y).
\end{equation*}
By linearity of $t$ with respect to its first argument (see \eqref{2.12}) and by \eqref{3.6.1}, we have
\begin{equation*}
|t(\xi+i\tau s\gamma,y)-t(\xi,y)|=|t(i\tau s\gamma,y)|\leq \lambda_3^{-1}s\tau,
\end{equation*}
which, together with \eqref{1.2} and \eqref{case1}, gives the estimate
\begin{eqnarray}\label{4.11}
|I_1|&\leq& \lambda^{-1}_3\lambda_0^{-1}s\tau\{|\partial_y\hat{v}|+\tau(\alpha_\pm+\beta |y|)|\hat{v}|+\lambda_3^{-1}(|\xi|+s\tau)|\hat{v}|\}\nonumber\\
 &\leq& C s\tau\{|\partial_y\hat{v}|+[\tau(\alpha_\pm+\beta |y|)+s\tau]|\hat{v}|\},
\end{eqnarray}
where $C$ depends on $\lambda_0$ only. On the other hand, by linearity of $t$ and by \eqref{3.7.1}, we obtain
\begin{equation*}
	\left|\partial_y\left(t(\xi+i\tau s\gamma,y)-t(\xi,y)\right)\right|=\left|\partial_y\left(t(i\tau s\gamma,y)\right)\right|\leq M_3s\tau,
\end{equation*}
which, together with \eqref{1.2}, \eqref{1.3} and \eqref{case1}, gives the estimate
\begin{equation}\label{4.12}
|I_2|\leq C s\tau\{|\partial_y\hat{v}|+[\tau(\alpha_\pm+\beta |y|)+s\tau]|\hat{v}|\},
\end{equation}
where $C$ depends on $\lambda_0$ and $M_0$ only.

Finally, by \eqref{2.14}, \eqref{3.1} and \eqref{case1},
\begin{equation}\label{4.8}
|I_3|=|2is\tau B_{\pm}(\xi,\gamma,y)-s^2\tau^2B_{\pm}(\gamma,\gamma,y)|
\leq C s\tau^2
\end{equation}
where $C$ depends only on $\lambda_0$. Putting together \eqref{4.7.1}, \eqref{4.12}, \eqref{4.11}, and \eqref{4.8} gives \eqref{4.7}.
\eproof
Lemma \ref{lem4.1} allows us to estimate $||P_{\pm}^0\hat{v}_{\pm}(\xi,\cdot)||_{L^2(\mathbb{R}_\pm)}$
instead of $||P_{\pm}\hat{v}_{\pm}(\xi,\cdot)||_{L^2(\mathbb{R}_\pm)}$. Let us now go further and note that,
similarly to \eqref{3.12} and \eqref{3.13}, we have
\begin{equation*}
P^0_+\hat{v}_+=E^0_{+}a_{nn}^+(y)F^0_+\hat{v}_+-\hat{v}_+\partial_y\big(a_{nn}^+(y)m_{+}(\xi,y)\big),
\end{equation*}
\begin{equation*}
P^0_-\hat{v}_-=F^0_{-}a_{nn}^+(y)E^0_-\hat{v}_-+\hat{v}_-\partial_y\big(a_{nn}^-(y)m_{-}(\xi,y)\big).
\end{equation*}
We can easily obtain, from \eqref{1.3} and \eqref{3.7}, that
\begin{equation}\label{4.17}
|P^0_+\hat{v}_+-E^0_{+}a_{nn}^+(y)F^0_+\hat{v}_+|\leq C |\xi||\hat{v}_+|
\end{equation}
and
\begin{equation}\label{4.18}
|P^0_-\hat{v}_--F^0_{-}a_{nn}^+(y)E^0_-\hat{v}_-|\leq C |\xi||\hat{v}_-|.
\end{equation}
where $C$ depends only on $\lambda_0$ and $M_0$.

\begin{lemma}\label{lem4.2}
Let $\tau\geq1$ and assume \eqref{case1}. There exists a positive constant $C$ depending only on $\lambda_0$ and $M_0$ such that,
if $s_0\leq 1/C$ then we have
\begin{eqnarray}\label{4.19}
&&\Lambda|a_{nn}^+(0)F^0_+\hat{v}_+(\xi,0)|^2+\Lambda^3|\hat{v}_+(\xi,0)|^2+\Lambda^4||\hat{v}_+(\xi,\cdot)||^2_{L^2(\mathbb{R}_+)}+\Lambda^2||\partial_y\hat{v}_+(\xi,\cdot)||^2_{L^2(\mathbb{R}_+)}\nonumber\\
&&\leq C ||P_{+}\hat{v}_{+}(\xi,\cdot)||^2_{L^2(\mathbb{R}_+)}
\end{eqnarray}
and
\begin{eqnarray}\label{4.20}
&&-\Lambda|a_{nn}^-(0)E^0_-\hat{v}_-(\xi,0)|^2-\Lambda^3|\hat{v}_-(\xi,0)|^2+\Lambda^4||\hat{v}_-(\xi,\cdot)||^2_{L^2(\mathbb{R}_-)}+\Lambda^2||\partial_y\hat{v}_-(\xi,\cdot)||^2_{L^2(\mathbb{R}_-)}\nonumber\\
&&\leq C ||P_{-}\hat{v}_{-}(\xi,\cdot)||^2_{L^2(\mathbb{R}_-)},
\end{eqnarray}
where ${\rm supp}(\hat{v}_+(\xi,\cdot))\subset[0,\frac{1}{\beta}]$ and  ${\rm supp}(\hat{v}_-(\xi,\cdot))\subset[-\frac{\alpha_-}{2\beta},0]$.
\end{lemma}
\pf. Since ${\rm supp}\, \hat{v}_+(x,y)$ is compact, $\hat{v}_+(\xi,y)\equiv 0$ when $|y|$ is large and the same holds for the function $\omega_+^0(\xi,y)$ defined in \eqref{omega0}. We now compute
\begin{eqnarray}\label{4.21}
&&||E^0_{+}\omega^0_{+}(\xi,\cdot)||^2_{L^2(\mathbb{R}_{+})}\nonumber\\
&&=\int_0^\infty|\partial_y\omega^0_+(\xi,y)+it_{+}(\xi,y)\omega^0_+(\xi,y)|^2dy+\int_0^\infty[\tau\alpha_{+}+\tau\beta y+m_{+}(\xi,y)]^2|\omega^0_+(\xi,y)|^2dy\nonumber\\
&&-2\Re \int_0^\infty[\tau\alpha_{+}+\tau\beta y+m_{+}(\xi,y)]\bar{\omega}^0_+(\xi,y)[\partial_y\omega^0_+(\xi,y)+it_{+}(\xi,y)\omega^0_+(\xi,y)]dy.
\end{eqnarray}
Integrating by parts, we easily get
\begin{equation}\label{4.22}
\begin{aligned}
&-2\Re \int_0^\infty[\tau\alpha_{+}+\tau\beta y+m_{+}(\xi,y)]\bar{\omega}^0_+(\xi,y)[\partial_y\omega^0_+(\xi,y)+it_{+}(\xi,y)\omega^0_+(\xi,y)]dy\\
&=[\tau\alpha_{+}+m_{+}(\xi,0)]|\omega^0_+(\xi,0)|^2+\int_0^\infty[\tau\beta +\partial_ym_{+}(\xi,y)]|\omega^0_+(\xi,y)|^2dy.
\end{aligned}
\end{equation}
By \eqref{case1} and \eqref{3.7}, we have that
\begin{equation}\label{4.23}
\tau\beta +\partial_ym_{+}(\xi,y)\geq\tau\beta -M_2|\xi|\geq\tau\beta -2\tau s_0\lambda_2^{-2}M_2\geq\tau\beta/2\geq 0
\end{equation}
provided $0<s_0\leq \frac{\beta\lambda_2^2}{4M_2}$. Combining \eqref{4.5}, \eqref{4.21}, \eqref{4.22} and \eqref{4.23} yields
\begin{eqnarray}\label{4.24}
||E^0_{+}\omega^0_{+}(\xi,\cdot)||^2_{L^2(\mathbb{R}_{+})}
&\geq&\int_0^\infty[\tau\alpha_{+}+\tau\beta y+m_{+}(\xi,y)]^2|\omega^0_+(\xi,y)|^2dy\nonumber\\&&+[\tau\alpha_{+}+m_{+}(\xi,0)]|\omega^0_+(\xi,0)|^2\nonumber
\\&\geq& C^{-1}\Lambda^2\int_0^\infty|\omega^0_+(\xi,y)|^2dy+C^{-1}\Lambda|\omega^0_+(\xi,0)|^2,
\end{eqnarray}
where $C$ depends only on $\lambda_0$.

Similarly, we have that
\begin{eqnarray}\label{4.25}
&&\lambda_0^{-2}||\omega^0_{+}(\xi,\cdot)||^2_{L^2(\mathbb{R}_{+})}
\geq\int_0^\infty|\partial_y\hat{v}_+(\xi,y)+it_{+}(\xi,y)\hat{v}_+(\xi,y)|^2dy\nonumber\\
&&+\int_0^\infty[\tau\alpha_{+}+\tau\beta y-m_{+}(\xi,y)]^2|\hat{v}_+(\xi,y)|^2dy
+[\tau\alpha_{+}-m_{+}(\xi,0)]|\hat{v}_+(\xi,0)|^2\nonumber\\
&&+\int_0^\infty[\tau\beta -\partial_ym_{+}(\xi,y)]|\hat{v}_+(\xi,y)|^2dy.
\end{eqnarray}
The assumption \eqref{case1} and \eqref{3.6} imply
\begin{equation*}
\tau\alpha_{+}+\tau\beta y-m_{+}(\xi,y)\geq\tau\alpha_+ -\lambda_2^{-1}|\xi|
\geq\tau\alpha_+ -2\lambda_2^{-3}\tau s_0
\geq\tau\alpha_+/2
\end{equation*}
provided $0<s_0\leq \frac{\alpha_+\lambda_2^3}{4}$. Thus, by choosing
\[
0<s_0\le\min\left\{1,\frac{\beta\lambda_2^2}{4M_2},\frac{\alpha_+\lambda_2^3}{4}\right\},
\]
we obtain from \eqref{4.23} and \eqref{4.25}
\begin{eqnarray}\label{4.27}
C||\omega^0_{+}(\xi,\cdot)||^2_{L^2(\mathbb{R}_{+})}
&\geq &\int_0^\infty|\partial_y\hat{v}_+(\xi,y)+it_{+}(\xi,y)\hat{v}_+(\xi,y)|^2dy\nonumber\\
&+&\Lambda^2\int_0^\infty|\hat{v}_+(\xi,y)|^2dy+\Lambda|\hat{v}_+(\xi,0)|^2.
\end{eqnarray}
Additionally, we can see that
\begin{eqnarray}\label{4.28}
&&\int_0^\infty|\partial_y\hat{v}_+(\xi,y)+it_{+}(\xi,y)\hat{v}_+(\xi,y)|^2dy\nonumber\\
&&\geq\eps\int_0^\infty\left(|\partial_y\hat{v}_+(\xi,y)|^2-2|\partial_y\hat{v}_+(\xi,y)||t_+(\xi,y)\hat{v}_+(\xi,y)|+|t_{+}(\xi,y)\hat{v}_+(\xi,y)|^2\right)dy\nonumber\\
&&\geq \eps\int_0^\infty\left(\frac 12|\partial_y\hat{v}_+(\xi,y)|^2-|t_+(\xi,y)|^2|\hat{v}_+(\xi,y)|^2\right)dy\nonumber\\
&&\geq\frac{\eps}{2}\int_0^\infty|\partial_y\hat{v}_+(\xi,y)|^2dy-\lambda_3^{-1}\eps|\xi|^2\int_0^\infty|\hat{v}_+(\xi,y)|^2dy,
\end{eqnarray}
for any $0<\eps<1$. Choosing $\eps$ sufficiently small, we obtain, from \eqref{4.27} and \eqref{4.28},
\begin{equation}\label{4.29}
C||\omega^0_{+}(\xi,\cdot)||^2_{L^2(\mathbb{R}_{+})}
\geq \int_0^\infty|\partial_y\hat{v}_+(\xi,y)|^2dy+\Lambda^2\int_0^\infty|\hat{v}_+(\xi,y)|^2dy+\Lambda|\hat{v}_+(\xi,0)|^2,
\end{equation}
where $C$ depends only on $\lambda_0$ and $M_0$.

Combining \eqref{4.24} and \eqref{4.29} yields
\begin{equation}\label{4.30}
\begin{aligned}
&\Lambda^2\int_0^\infty|\partial_y\hat{v}_+(\xi,y)|^2+\Lambda^4\int_0^\infty|\hat{v}_+(\xi,y)|^2+\Lambda^3|\hat{v}_+(\xi,0)|^2+\Lambda|\omega^0_+(\xi,0)|^2\\
\leq&C||E^0_{+}\omega^0_{+}(\xi,\cdot)||^2_{L^2(\mathbb{R}_{+})},
\end{aligned}
\end{equation}
where $C$ depends only on $\lambda_0$ and $M_0$. From \eqref{4.7}, since ${\rm supp}(\hat{v}_+(\xi,\cdot))\subset [0,1/\beta]$ and \eqref{case1} holds, we have
\begin{eqnarray}\label{N*}
	||P_{+}^0\hat{v}_{+}(\xi,\cdot)||^2_{L^2(\mathbb{R}_+)}&\leq& 2 ||P_{+}^0\hat{v}_{+}(\xi,\cdot)||_{L^2(\mathbb{R}_+)}\nonumber\\&&+
	Cs_0^2\left(\Lambda^2\int_0^\infty|\partial_y\hat{v}_+(\xi,y)|^2+\Lambda^4\int_0^\infty|\hat{v}_+(\xi,y)|^2\right).
\end{eqnarray}
Moreover, by \eqref{4.17} and \eqref{case1},
\begin{equation}\label{N1*}
	||E^0_{+}\omega^0_{+}(\xi,\cdot)||^2_{L^2(\mathbb{R}_{+})}\leq 2||P_{+}^0\hat{v}_{+}(\xi,\cdot)||^2_{L^2(\mathbb{R}_+)}+
	Cs_0^2\Lambda^2\int_0^\infty|\hat{v}_+(\xi,y)|^2.
\end{equation}
Finally, by \eqref{4.30}, \eqref{N*} and \eqref{N1*} we get \eqref{4.19}, provided $s_0$ is small enough.

Now, we proceed to prove \eqref{4.20}. Applying the same arguments leading to \eqref{4.22}, we have that
\begin{equation}\label{4.31}
\begin{aligned}
&||F^0_{-}\omega^0_{-}(\xi,\cdot)||^2_{L^2(\mathbb{R}_{-})}\\
\geq&\int^0_{-\infty}[\tau\alpha_{-}+\tau\beta y-m_{-}(\xi,y)]^2|\omega^0_-(\xi,y)|^2dy-[\tau\alpha_{-}-m_{-}(\xi,0)]|\omega^0_-(\xi,0)|^2\\
&+\int^0_{-\infty}[\tau\beta-\partial_ym_{-}(\xi;y)]|\omega^0_-(\xi,y)|^2dy.
\end{aligned}
\end{equation}
By \eqref{3.6} and \eqref{case1} and since ${\rm supp}(\hat{v}_-(\xi,\cdot))\subset[-\frac{\alpha_-}{2\beta},0]$,  we can see that
\begin{equation}\label{4.32}
\tau\alpha_{-}+\tau\beta y-m_{-}(\xi,y)\geq\tau\alpha_-/2 -\lambda_2^{-1}|\xi|\geq\tau\alpha_-/2 -2\lambda_2^{-3}\tau s_0\geq\tau\alpha_-/4
\end{equation}
provided $0<s_0\leq \frac{\alpha_-\lambda_2^3}{8}$. From \eqref{4.31} and \eqref{4.32}, it follows
\begin{eqnarray}\label{4.33}
||F^0_{-}\omega^0_{-}(\xi,\cdot)||^2_{L^2(\mathbb{R}_{-})}
&\geq&\frac{\alpha_-^2}{16}\tau^2\int_0^\infty|\omega^0_-(\xi,y)|^2dy-\tau\alpha_-|\omega^0_-(\xi,0)|^2\nonumber\\&\geq& C\Lambda^2\int_0^\infty|\omega^0_-(\xi,y)|^2dy-C\Lambda|\omega^0_-(\xi,0)|^2.
\end{eqnarray}
Arguing as before and recalling \eqref{4.5} we obtain \eqref{4.20}.
\eproof

We now take into account the transmission conditions. 
\begin{lemma}\label{lem4.3}
Let  $\tau\geq1$ and assume \eqref{case1}. There exists a positive constant $C$ depending only on $\lambda_0$ and $M_0$ such that if $s_0\leq 1/C$ then
\begin{eqnarray}\label{4.37}
\Lambda\sum_{\pm}|V_{\pm}(\xi)|^2\!\!\!\!&+\!\!\!\!&\Lambda^3\sum_{\pm}|\hat{v}_{\pm}(\xi,0)|^2+\Lambda^4\sum_{\pm}||\hat{v}_{\pm}(\xi,\cdot)||^2_{L^2(\mathbb{R}_{\pm})}
+\Lambda^2\sum_{\pm}||\partial_y\hat{v}_{\pm}(\xi,\cdot)||^2_{L^2(\mathbb{R}_{\pm})}\nonumber\\
&&\leq C\sum_{\pm}||P_{\pm}\hat{v}_{\pm}(\xi,\cdot)||^2_{L^2(\mathbb{R}_{\pm})}+C\Lambda|\eta_1(\xi)|^2+C\Lambda^3|\eta_0(\xi)|^2,
\end{eqnarray}
where ${\rm supp}(\hat{v}_\pm(\xi,\cdot))\subset[-\frac{c_0}{2\beta},\frac{c_0}{\beta}]$ with  $c_0=\min{(\alpha_-,1)}$. \end{lemma}
\pf. It follows from \eqref{4.19} and \eqref{omega0} that, for some $C$ depending only on $\lambda_0$ and $M_0$,
\begin{equation}\label{4.38}
\Lambda|\omega_0^+(\xi,0)|^2+\Lambda^3|\hat{v}_+(\xi,0)|^2\leq C ||P_{+}\hat{v}_{+}(\xi,\cdot)||^2_{L^2(\mathbb{R}_{+})}.
\end{equation}
By \eqref{2.20}, \eqref{omega0}, \eqref{3.6} and \eqref{3.6.1} we easily get
\[V_+(\xi)=\omega_0^+(\xi,0)-a_{nn}^+(0)(\tau s t_+(\gamma,0)+m_+(\xi,0))\hat{v}_+(\xi,0),\]
and hence
\begin{equation}\label{4.39}
\Lambda|V_+(\xi)|^2\leq2\Lambda|\omega_0^+(\xi,0)|^2+C\Lambda^3|\hat{v}_+(\xi,0)|^2\leq C||P_{+}\hat{v}_{+}(\xi,\cdot)||^2_{L^2(\mathbb{R}_{+})},
\end{equation}
where $C$ depends only on $\lambda_0$ and $M_0$.

By \eqref{2.18} and \eqref{4.19}, we have that
\begin{equation}\label{4.40}
\Lambda^3|\hat{v}_-(\xi,0)|^2\leq 2\Lambda^3|\hat{v}_+(\xi,0)|^2+2\Lambda^3|\eta_0(\xi)|^2
\leq C ||P_{+}\hat{v}_{+}(\xi,\cdot)||^2_{L^2(\mathbb{R}_{+})}+2\Lambda^3|\eta_0(\xi)|^2.
\end{equation}
Using the definition of $\eta_1$ (see \eqref{2.19}) and \eqref{4.39}, we also deduce that
\begin{equation}\label{4.42}
\Lambda|V_-(\xi)|^2\leq 2\Lambda|V_+(\xi)|^2+2\Lambda|\eta_1(\xi)|^2\leq C ||P_{+}\hat{v}_{+}(\xi,\cdot)||^2_{L^2(\mathbb{R}_{+})}+2\Lambda|\eta_1(\xi)|^2.
\end{equation}
Putting together \eqref{4.38}, \eqref{4.39}, \eqref{4.40} and \eqref{4.42} implies
\begin{equation}\label{4.43}
\Lambda^3\sum_{\pm}|\hat{v}_{\pm}(\xi,0)|^2+\Lambda\sum_{\pm}|V_{\pm}(\xi)|^2\leq C||P_{+}\hat{v}_{+}(\xi,0)||^2_{L^2(\mathbb{R}_{+})}+2\Lambda^3|\eta_0(\xi)|^2+2\Lambda|\eta_1(\xi)|^2.
\end{equation}
We now use \eqref{4.19} and  \eqref{4.20} and get
\begin{equation*}
\begin{aligned}
&\Lambda^4\sum_{\pm}||\hat{v}_{\pm}(\xi,\cdot)||_{L^2(\mathbb{R}_{\pm})}+\Lambda^2\sum_{\pm}||\partial_y\hat{v}_{\pm}(\xi,\cdot)||_{L^2(\mathbb{R}_{\pm})}\\
\leq& C\sum_{\pm}||P_{\pm}\hat{v}_{\pm}(\xi,\cdot)||^2_{L^2(\mathbb{R}_{\pm})}+\Lambda|\omega_0^-(\xi,0)|^2+\Lambda^3|\hat{v}_-(\xi,0)|^2
\end{aligned}
\end{equation*}
Arguing similarly as we did for \eqref{4.39} and using \eqref{4.42} and \eqref{4.43}, we get
\begin{equation}\label{4.44}
\begin{aligned}
&\Lambda^4\sum_{\pm}||\hat{v}_{\pm}(\xi,\cdot)||_{L^2(\mathbb{R}_{\pm})}+\Lambda^2\sum_{\pm}||\partial_y\hat{v}_{\pm}(\xi,\cdot)||_{L^2(\mathbb{R}_{\pm})}\\
\leq & C\sum_{\pm}||P_{\pm}\hat{v}_{\pm}(\xi,\cdot)||^2_{L^2(\mathbb{R}_{\pm})}+2\Lambda|V_-(\xi)|^2+C\Lambda^3|\hat{v}_-(\xi,0)|^2\\
\leq&C\left( \sum_{\pm}||P_{\pm}\hat{v}_{\pm}(\xi,\cdot)||^2_{L^2(\mathbb{R}_{\pm})}+\Lambda|\eta_1(\xi)|^2+\Lambda^3|\eta_0(\xi)|^2\right),
\end{aligned}
\end{equation}
where $C$ depends on $\lambda_0$ and $M_0$ only.
The proof is complete by combining \eqref{4.43} and \eqref{4.44}.
\eproof

Since $\tau\ge 1$, it is easily seen that \eqref{4.37} implies
\begin{equation}\label{4.377}
\begin{aligned}
&\Lambda\sum_{\pm}|V_{\pm}(\xi)|^2+\Lambda^3\sum_{\pm}|\hat{v}_{\pm}(\xi,0)|^2+\frac{\Lambda^4}{\tau}\sum_{\pm}||\hat{v}_{\pm}(\xi,\cdot)||^2_{L^2(\mathbb{R}_{\pm})}
+\frac{\Lambda^2}{\tau}\sum_{\pm}||\partial_y\hat{v}_{\pm}(\xi,\cdot)||^2_{L^2(\mathbb{R}_{\pm})}\\
\leq& C\left(\sum_{\pm}||P_{\pm}\hat{v}_{\pm}(\xi,\cdot)||^2_{L^2(\mathbb{R}_{\pm})}+\Lambda|\eta_1(\xi)|^2+\Lambda^3|\eta_0(\xi)|^2\right),
\end{aligned}
\end{equation}
where $C$ depends on $\lambda_0$ and $M_0$ only.

\bigskip
\noindent{\bf Case 2:}
\begin{equation}\label{case2}
	\frac{m_+(\xi,0)}{(1-\kappa)\,\alpha_+}\leq\tau\leq\frac{\lambda_2^2|\xi|}{2s_0}.
\end{equation}

\bigskip
In this case, \eqref{3.6} implies
\begin{equation}\label{4.466}
\frac{\lambda_2|\xi|}{\alpha_+}\le\tau\le\frac{\lambda_2^2|\xi|}{2s_0}.
\end{equation}
In addition, in view of the definition of $\zeta_\pm$, \eqref{3.1}, \eqref{case2}, and recalling that
$\lambda_2=\sqrt{\lambda_0\lambda_1}$ and $s\leq s_0$, we have that
\begin{equation}\label{4.47}
\Re\zeta_{\pm}
\geq \frac{3}{4}\lambda_2^2|\xi|^2.
\end{equation}
It is not hard to see from \eqref{3.8}, \eqref{3.9}, \eqref{4.466}, \eqref{4.47} that
\begin{equation}\label{4.49}
|\partial_y\sqrt{\zeta_{\pm}}|\leq M_5|\xi|,
\end{equation}
where $M_5$ depends only on $\lambda_0$ and $M_0$. Moreover, if we set $R_{\pm}=\Re\sqrt{\zeta_{\pm}}\geq 0$ and
$J_{\pm}=\Im\sqrt{\zeta_{\pm}}$, then \eqref{4.49} gives
\begin{equation}\label{4.50}
|\partial_yR_{\pm}|+|\partial_yJ_{\pm}|\leq M_5|\xi|.
\end{equation}
Using \eqref{4.49}, we can easily obtain from \eqref{3.12}, \eqref{3.13} that
\begin{equation}\label{4.51}
|P_+\hat{v}_+(\xi,y)-E_{+}a_{nn}^+(y)F_+\hat{v}_+(\xi,y)|\leq C|\xi||\hat{v}_+(\xi,y)|
\end{equation}
and
\begin{equation}\label{4.52}
|P_-\hat{v}_-(\xi,y)-F_{-}a_{nn}^-(y)E_-\hat{v}_-(\xi,y)|\leq C|\xi||\hat{v}_-(\xi,y)|,
\end{equation}
where $C$ depends only on $\lambda_0$ and $M_0$.

We now prove the following lemma.
\begin{lemma}\label{lem4.4}
Assume \eqref{case2}. There exists a  positive constant $C$ depending only on $\lambda_0$ and $M_0$ such that,
if $0<s_0\leq C^{-1}$, $\beta\geq C$ and $\tau\geq C$, then we have
\begin{equation}\label{4.53}
\begin{aligned}
&\Lambda|V_+(\xi)+a_{nn}^+(0)\sqrt{\zeta_{+}(\xi,0)}\hat{v}_+(\xi,0)|^2+\Lambda^2||a_{nn}^+(y)F_+\hat{v}_+(\xi,\cdot)||^2_{L^2(\mathbb{R}_{+})}\\
\leq&C||E_{+}a_{nn}^+(y)F_+\hat{v}_+(\xi,\cdot)||^2_{L^2(\mathbb{R}_{+})}
\end{aligned}
\end{equation}
and
\begin{equation}\label{4.54}
\begin{aligned}
&\Lambda|V_+(\xi)+a_{nn}^+(0)\sqrt{\zeta_{+}(\xi,0)}\hat{v}_+(\xi,0)|^2+\Lambda^3|\hat{v}_+(\xi,0)|^2\\
&+\Lambda^4\int_0^\infty|\hat{v}_+(\xi,y)|^2dy+\Lambda^2\int_0^\infty|\partial_y\hat{v}_+(\xi,y)|^2dy\leq C||P_{+}\hat{v}_{+}(\xi,\cdot)||^2_{L^2(\mathbb{R}_{+})}
\end{aligned}
\end{equation}
provided ${\rm supp}(\hat{v}_+(\xi,\cdot))\subset[0,\frac{1}{\beta}]$.
\end{lemma}
\pf. We write
\[
E_{+}\omega_{+}(\xi,y)=[\partial_y+it_{+}(\xi+i\tau s\gamma,y)-\tau\varphi'_{+}-\sqrt{\zeta_{+}}]\omega_+(\xi,y):=I_3-I_4,
\]
where $I_3=\partial_y\omega_++it_{+}(\xi+i\tau s\gamma,y)\omega_+-iJ_+\omega_+$ and
$I_4=\tau\alpha_{+}\omega_++\tau\beta y\omega_++R_+\omega_+$. Our task now is to estimate
\begin{equation}\label{4.55}
||E_{+}\omega_{+}(\xi,\cdot)||^2_{L^2(\mathbb{R}_{+})}
=\int_0^\infty|I_3|^2dy+\int_0^\infty[\tau\alpha_{+}+\tau\beta y+R_{+}]^2|\omega_+|^2dy-2\Re \int_0^\infty I_3\bar{I_4}dy.
\end{equation}
Observe that
\begin{equation}\label{4.56}
\begin{aligned}
-2\Re \int_0^\infty I_3\bar{I_4}=&-\int_0^\infty[\tau\alpha_{+}+\tau\beta y+R_{+}(\xi,y)]\partial_y(|\omega_+(\xi,y)|^2)dy\\
&+2\int_0^\infty[\tau\alpha_{+}+\tau\beta y+R_{+}(\xi,y)]t_{+}(\tau s\gamma,y)|\omega_+(\xi,y)|^2dy\\
=&\int_0^\infty[\tau\beta +\partial_yR_{+}(\xi,y)]|\omega_+(\xi,y)|^2dy+[\tau\alpha_{+}+R_{+}(\xi,0)]|\omega_+(\xi,0)|^2\\
&+2\int_0^\infty[\tau\alpha_{+}+\tau\beta y+R_{+}(\xi,y)]t_{+}(\tau s\gamma,y)|\omega_+(\xi,y)|^2dy\\
\geq&\int_0^\infty[\tau\beta +\partial_yR_{+}(\xi,y)-\lambda_3^{-1}s\tau(\tau\alpha_++\tau\beta y+R_+)]|\omega_+(\xi,y)|^2dy\\
&+[\tau\alpha_{+}+R_{+}(\xi,0)]|\omega_+(\xi,0)|^2,
\end{aligned}
\end{equation}
where in the last inequality we have used the fact that $R_+\ge 0$. Combining \eqref{4.55} and \eqref{4.56} yields
\begin{equation}\label{4.57}
\begin{aligned}
&||E_{+}\omega_{+}(\xi,\cdot)||^2_{L^2(\mathbb{R}_{+})}\\
\geq&\int_0^\infty[(\tau\alpha_{+}+\tau\beta y+R_{+})^2+\tau\beta +\partial_yR_{+}(\xi,y)-\lambda_3^{-1}s\tau(\tau\alpha_{+}+\tau\beta y+R_{+})]|\omega_+(\xi,y)|^2dy\\
&+[\tau\alpha_{+}+R_{+}(\xi,0)]|\omega_+(\xi,0)|^2\\
\geq& \frac{\Lambda^2}{C}\int_0^\infty|\omega_+(\xi,y)|^2dy+\frac{\Lambda}{C}|\omega_+(\xi,0)|^2
\end{aligned}
\end{equation}
provided $s_0$ is small enough. Formulas \eqref{2.20} and \eqref{3.11} give
\begin{equation}\label{4.58}
\omega_+(\xi,0)=V_+(\xi)+a_{nn}^+(0)\sqrt{\zeta_{+}(\xi,0)}\hat{v}_+(\xi,0),
\end{equation}
which leads to \eqref{4.53}  by \eqref{4.57}.

We now want to derive \eqref{4.54}. Let us write
$$F_{+}\hat{v}_{+}=[\partial_y+it_{+}(\xi+i\tau s\gamma;y)-\tau\varphi'_{+}+\sqrt{\zeta_{+}}]\hat{v}_+:=I_5-I_6,$$
where $I_5=\partial_y\hat{v}_++it_{+}(\xi+i\tau s\gamma;y)\hat{v}_++iJ_+\hat{v}_+$ and
$I_6=\tau\alpha_{+}\hat{v}_++\tau\beta y\hat{v}_+-R_+\hat{v}_+$. Thus, we have
\begin{equation}\label{4.59}
\begin{aligned}
&||F_{+}\hat{v}_{+}(\xi,\cdot)||^2_{L^2(\mathbb{R}_{+})}\\
=&\int_0^\infty|I_5|^2dy+\int_0^\infty[\tau\alpha_{+}+\tau\beta y-R_{+}]^2|\hat{v}_+(\xi,y)|^2dy-2\Re \int_0^\infty I_5\bar{I_6}dy.
\end{aligned}
\end{equation}
Repeating the computations of \eqref{4.56} and \eqref{4.57} yields
\begin{equation}\label{4.60}
\begin{aligned}
&||F_{+}\hat{v}_{+}(\xi,\cdot)||^2_{L^2(\mathbb{R}_{+})}\\
\geq&\int_0^\infty|I_5|^2dy+\int_0^\infty[\tau\alpha_{+}+\tau\beta y-R_{+}]^2|\hat{v}_+(\xi,y)|^2dy+\int_0^\infty(\tau\beta -\partial_yR_{+})|\hat{v}_+(\xi,y)|^2dy\\
&-Cs\tau\int_0^\infty|\tau\alpha_{+}+\tau\beta y-R_{+}||\hat{v}_+(\xi,y)|^2dy+[\tau\alpha_{+}-R_{+}(\xi,0)]|\hat{v}_+(\xi,0)|^2.
\end{aligned}
\end{equation}
We observe that
\[R_+^2=\frac{\Re\zeta_++|\zeta_+|}{2}\]
and, by simple calculations,
\begin{equation}\label{4.61}
|\zeta_{\pm}|
\leq -\Re \zeta_{\pm}+2\frac{B_{\pm}(\xi,\xi,y)}{a_{nn}^{\pm}(y)},
\end{equation}
which gives the estimate
\begin{equation}\label{4.62}
R_+(\xi,y)\leq\sqrt{\frac{B_{+}(\xi,\xi;y)}{a_{nn}^{+}(y)}}=m_+(\xi,y).
\end{equation}
From \eqref{case2} and \eqref{4.62}, we deduce that
\begin{equation}\label{4.63}
\tau\alpha_{+}-R_{+}(\xi,0)\geq\tau\alpha_{+}-m_{+}(\xi,0)\geq\tau\alpha_{+}-(1-\kappa)\tau\alpha_{+}=\kappa\tau\alpha_{+}.
\end{equation}
On the other hand, using \eqref{4.63}, \eqref{4.50} and \eqref{4.466}, we can obtain that for $y\ge 0$
\begin{equation*}
\begin{aligned}
\tau\alpha_{+}+\tau\beta y-R_{+}(\xi,y)=&\tau\alpha_{+}-R_{+}(\xi,0)+\tau\beta y-R_{+}(\xi,y)+R_{+}(\xi,0)\\
\geq&\kappa\tau\alpha_{+}+y(\tau\beta-C\tau)
\geq\kappa\tau\alpha_{+}
\end{aligned}
\end{equation*}
provided $\beta$ is large enough. Furthermore, if $0\leq y\leq 1/\beta$, then
\begin{equation}\label{4.65}
[\tau\alpha_{+}+\tau\beta y-R_{+}]^2+(\tau\beta -\partial_yR_{+})-Cs\tau|\tau\alpha_{+}+\tau\beta y-R_{+}|
\geq (\kappa\tau\alpha_{+})^2/4
\end{equation}
provided $s_0$ is small enough and $\tau$ is large enough. Now it follows from \eqref{4.60}, \eqref{4.63}, and \eqref{4.65}
and arguing as in \eqref{4.28}, that
\begin{equation}\label{4.66}
\begin{aligned}
&C||F_{+}\hat{v}_{+}(\xi,y)||^2_{L^2(\mathbb{R}_{+})}\\
\geq&\int_0^\infty|\partial_y\hat{v}_+(\xi,y)|^2dy+\Lambda^2\int_0^\infty|\hat{v}_+(\xi,y)|^2dy+\Lambda|\hat{v}_+(\xi,0)|^2.
\end{aligned}
\end{equation}
Finally, by \eqref{4.51}, \eqref{4.53}, and \eqref{4.66}, we can easily derive \eqref{4.54} provided $\beta\geq C$, $\tau\ge C$ and $s_0\leq 1/C$ for some $C$ depending on $\lambda_0$ and $M_0$.
\eproof

Similarly, we can prove that
\begin{lemma}\label{lem4.5}
Assume \eqref{case2}. There exists a positive constant $C$ depending only on $\lambda_0$ and $M_0$ such that,
if $0<s_0\leq C^{-1}$ and $\tau\geq C$ then we have
\begin{equation}\label{4.67}
\begin{aligned}
&-\Lambda|V_-(\xi)-a_{nn}^-(0)\sqrt{\zeta_{-}}\hat{v}_-(\xi,0)|^2+\Lambda||a_{nn}^-(y)E_-\hat{v}_-(\xi,\cdot)||^2_{L^2(\mathbb{R}_{-})}\\
\leq &C||F_{-}a_{nn}^-(y)E_-\hat{v}_-(\xi,\cdot)||^2_{L^2(\mathbb{R}_{-})}
\end{aligned}
\end{equation}
and
\begin{equation}\label{4.68}
\begin{aligned}
&-\Lambda|V_-(\xi)-a_{nn}^-(0)\sqrt{\zeta_{-}}\hat{v}_-(\xi,0)|^2-\Lambda^3|\hat{v}_-(\xi,0)|^2+\Lambda^3\int_{-\infty}^0|\hat{v}_-(\xi,y)|^2dy\\
&+\Lambda\int_{-\infty}^0|\partial_y\hat{v}_-(\xi,y)|^2dy\leq C ||P_{-}\hat{v}_{-}(\xi,\cdot)||^2_{L^2(\mathbb{R}_{-})},
\end{aligned}
\end{equation}
provided ${\rm supp}(\hat{v}_-(\xi,\cdot))\subset[-\frac{\alpha_-}{2\beta},0]$.
\end{lemma}
\pf. Let $\omega_-(\xi,y)=a_{nn}^-(y)E_-\hat{v}_-(\xi,y)=a_{nn}^-(y)[\partial_y+it_{-}(\xi+i\tau s\gamma,y)-\tau\varphi'_{-}-\sqrt{\zeta_{-}}]\hat{v}_-(\xi,y)$.
If we write
\[F_{-}\omega_{-}(\xi,y)=
I_7-I_8,\]
where
\begin{eqnarray*}
I_7&=&\partial_y\omega_-+it_{-}(\xi,y)\omega_-+iJ_-\omega_-\\
I_8&=&\tau\alpha_{-}\omega_-+\tau\beta y\omega_-+t_{-}(\tau s\gamma,y)\omega_--R_-\omega_-,
\end{eqnarray*}
we have
\begin{equation*}
\begin{aligned}
&||F_{-}\omega_{-}(\xi,\cdot)||^2_{L^2(\mathbb{R}_{-})}\geq-2\Re \int_{-\infty}^0 I_7\bar{I_8}dy\\
=&-\int_{-\infty}^0[\tau\alpha_{-}+\tau\beta y+t_{-}(\tau s\gamma,y)-R_-(\xi,y)]\partial_y(|\omega_-(\xi,y)|^2)dy\\
=&\int_{-\infty}^0[\tau\beta+\partial_y t_{-}(\tau s\gamma,y)-\partial_yR_-(\xi,y)]|\omega_-(\xi,y)|^2dy\\
&-[t_{-}(\tau s\gamma,0)+\tau\alpha_{-}-R_-(\xi,0)]|\omega_-(\xi,0)|^2\\
\geq&\int_{-\infty}^0\tau[\beta-M_3s -2M_5s_0\lambda_2^{-2}]|\omega_-(\xi,y)|^2dy-(\lambda_3 s+\alpha_+)\tau|\omega_-(\xi,0)|^2,
\end{aligned}
\end{equation*}
hence, by  \eqref{4.466},
\begin{equation}\label{4.69}
||F_{-}\omega_{-}(\xi,\cdot)||^2_{L^2(\mathbb{R}_{-})}\geq C\Lambda\int_{-\infty}^0|\omega_-(\xi,y)|^2dy-C\Lambda|\omega_-(\xi,0)|^2,
\end{equation}
provided $s_0$ is small enough. Since, by \eqref{2.20} and \eqref{3.10},
\begin{equation*}
\omega_-(\xi,0)=V_-(\xi)-a_{nn}^-(0)\sqrt{\zeta_{-}}\hat{v}_-(\xi,0),
\end{equation*}
we get \eqref{4.67}.

To derive \eqref{4.68}, we denote
$$E_-\hat{v}_-(\xi,y)=I_9-I_{10},$$
where
\begin{eqnarray*}I_9&=&\partial_y\hat{v}_-+it_{-}(\xi,y)\hat{v}_--iJ_-\hat{v}_-,\\
I_{10}&=&\tau\alpha_{-}\hat{v}_-+\tau\beta y\hat{v}_-+t_{-}(\tau s\gamma,y)\hat{v}_-+R_-\hat{v}_-.
\end{eqnarray*}
 Observe that if $-\frac{\alpha_-}{2\beta}\leq y\leq 0$ then
\begin{equation}\label{4.555}
\tau\alpha_{-}+\tau\beta y+t_{-}(\tau s\gamma,y)+R_-
\geq\tau\alpha_{-}/2-\lambda_3^{-1}s\tau\geq \tau\alpha_{-}/4
\end{equation}
provided $s_0$ is small. Furthermore, by choosing again $s_0$ small, we can make
\begin{equation}\label{4.556}
\tau\beta+\partial_yR_-+\partial_yt_{-}(\tau s\gamma,y)\ge \tau\left(\beta -2M_5s_0\lambda_2^{-2}-M_3s_0\right)\ge 0.
\end{equation}
With the help of \eqref{4.555} and \eqref{4.556}, and arguing as in \eqref{4.28} we get
\begin{equation}\label{4.71}
\begin{aligned}
&C||E_-\hat{v}_-(\xi,\cdot)||^2_{L^2(\mathbb{R}_{-})}\\
\geq&\int_{-\infty}^0|\partial_y\hat{v}_-(\xi,y)|^2dy+\Lambda^2\int_{-\infty}^0|\hat{v}_-(\xi,y)|^2dy-\Lambda|\hat{v}_-(\xi,0)|^2.
\end{aligned}
\end{equation}
Using  \eqref{4.67}, \eqref{4.71} and \eqref{4.52}, we obtain \eqref{4.68} provided $\tau$ is large.
\eproof

\begin{lemma}\label{lem4.6}
Assume \eqref{case2}. There exists a positive constant $C$,  depending only on $\lambda_0$ and $M_0$, such that if $s_0\leq C^{-1}$, $\beta\geq C$ and $\tau\ge C$ then we have
\begin{equation}\label{4.72}
\begin{aligned}
&\Lambda\sum_{\pm}|V_{\pm}(\xi)|^2+\Lambda^3\sum_{\pm}|\hat{v}_{\pm}(\xi,0)|^2+\Lambda^3\sum_{\pm}||\hat{v}_{\pm}(\xi,\cdot)||^2_{L^2(\mathbb{R}_{\pm})}
+\Lambda\sum_{\pm}||\partial_y\hat{v}_{\pm}(\xi,\cdot)||^2_{L^2(\mathbb{R}_{\pm})}\\
\leq &C\left(\sum_{\pm}||P_{\pm}\hat{v}_{\pm}(\xi,\cdot)||^2_{L^2(\mathbb{R}_{\pm})}+\Lambda|\eta_1(\xi)|^2+\Lambda^3|\eta_0(\xi)|^2\right),
\end{aligned}
\end{equation}
provided ${\rm supp}(\hat{v}_\pm(\xi,\cdot))\subset[-\frac{c_0}{2\beta},\frac{c_0}{\beta}]$ with  $c_0=\min{(\alpha_-,1)}$.
\end{lemma}
\pf. We obtain from \eqref{4.54} that
\begin{equation}\label{4.73}
\Lambda|\omega_+(\xi,0)|^2+\Lambda^3|\hat{v}_+(\xi,0)|^2\leq C ||P_{+}\hat{v}_{+}(\xi,\cdot)||^2_{L^2(\mathbb{R}_{+})}.
\end{equation}
On the other hand,
\begin{equation}\label{4.74}
\Lambda|V_+(\xi)|^2\le 2\Lambda|\omega_+(\xi,0)|^2+C\Lambda^3|\hat{v}_+(\xi,0)|^2\leq C||P_{+}\hat{v}_{+}(\xi,\cdot)||^2_{L^2(\mathbb{R}_{+})}.
\end{equation}
Using the definition of $\eta_0$ and \eqref{4.73}, we see that
\begin{equation}\label{4.75}
\Lambda^3|\hat{v}_-(\xi,0)|^2\leq 2\Lambda^3|\hat{v}_+(\xi,0)|^2+ 2\Lambda^3|\eta_0(\xi)|^2\leq C||P_{+}\hat{v}_{+}(\xi,\cdot)||^2_{L^2(\mathbb{R}_{+})}+2\Lambda^3|\eta_0(\xi)|^2.
\end{equation}
Summing up \eqref{4.73} and \eqref{4.75} yields
\begin{equation}\label{4.76}
\Lambda^3\sum_{\pm}|\hat{v}_{\pm}(\xi,0)|^2\leq C||P_{+}\hat{v}_{+}(\xi,\cdot)||^2_{L^2(\mathbb{R}_{+})}+2\Lambda^3|\eta_0(\xi)|^2.
\end{equation}
Likewise, the definition of $\eta_1$ and \eqref{4.74} lead to
\begin{equation}\label{4.77}
\Lambda|V_-(\xi)|^2\leq C ||P_{+}\hat{v}_{+}(\xi,\cdot)||^2_{L^2(\mathbb{R}_{+})}+2\Lambda|\eta_1(\xi)|^2.
\end{equation}
Putting together \eqref{4.74}, \eqref{4.76}, and \eqref{4.77}, we deduce that
\begin{equation}\label{4.78}
\Lambda^3\sum_{\pm}|\hat{v}_{\pm}(\xi,0)|^2+\Lambda\sum_{\pm}|V_{\pm}(\xi)|^2\leq C||P_{+}\hat{v}_{+}(\xi,\cdot)||^2_{L^2(\mathbb{R}_{+})}+2\Lambda^3|\eta_0(\xi)|^2+2\Lambda|\eta_1(\xi)|^2.
\end{equation}
Finally, we first use \eqref{4.54}, recall that $\Lambda\geq\tau\geq1$, \eqref{4.68}, and then \eqref{4.77}, \eqref{4.78} to derive
\begin{equation}\label{4.79}
\begin{aligned}
&\Lambda^3\sum_{\pm}||\hat{v}_{\pm}(\xi,\cdot)||^2_{L^2(\mathbb{R}_{\pm})}+\Lambda\sum_{\pm}||\partial_y\hat{v}_{\pm}(\xi,\cdot)||^2_{L^2(\mathbb{R}_{\pm})}\\
\leq &C\sum_{\pm}||P_{\pm}\hat{v}_{\pm}(\xi,\cdot)||^2_{L^2(\mathbb{R}_{\pm})}+\Lambda|\hat{V}_-(\xi)-a_{nn}^-(0)R_-(\xi,0)\hat{v}_-(\xi,0)|^2+\Lambda^3|\hat{v}_-(\xi,0)|^2\\
\leq&C\left( \sum_{\pm}||P_{\pm}\hat{v}_{\pm}(\xi,\cdot)||^2_{L^2(\mathbb{R}_{\pm})}+\Lambda|\eta_1(\xi)|^2+\Lambda^3|\eta_0(\xi)|^2\right).
\end{aligned}
\end{equation}
The proof is complete by combining \eqref{4.78} and \eqref{4.79}.
\eproof

We conclude Case 2 by observing that \eqref{4.72} can be written in the form
\begin{equation}\label{4.80}
\begin{aligned}
&\Lambda\sum_{\pm}|V_{\pm}(\xi)|^2+\Lambda^3\sum_{\pm}|\hat{v}_{\pm}(\xi,0)|^2+\frac{1}{\tau}\big(\Lambda^4\sum_{\pm}||\hat{v}_{\pm}(\xi,\cdot)||^2_{L^2(\mathbb{R}_{\pm})}
+\Lambda^2\sum_{\pm}||\partial_y\hat{v}_{\pm}(\xi,\cdot)||^2_{L^2(\mathbb{R}_{\pm})}\big)\\
\leq&C\left( \sum_{\pm}||P_{\pm}\hat{v}_{\pm}(\xi,\cdot)||^2_{L^2(\mathbb{R}_{\pm})}+\Lambda|\eta_1(\xi)|^2+\Lambda^3|\eta_0(\xi)|^2\right),
\end{aligned}
\end{equation}
where $C$ depends only on $\lambda_0$ and $M_0$.

\bigskip
\noindent{\bf Case 3:}
\begin{equation}\label{case3}
	\tau\leq\frac{m_+(\xi,0)}{(1-\kappa)\,\alpha_+}.
\end{equation}

In this case, we have
\[\tau\leq \frac{2\lambda_2^{-1}|\xi|}{\alpha_++L\alpha_-}\quad\text{(from \eqref{3.6}, \eqref{kappa})}. \]
From the definition of $\zeta_\pm$ (see \eqref{2.14}) and the inequality
\begin{equation*}
B_{\pm}(\xi,\xi;y)-s^2\tau^2B_{\pm}(\gamma,\gamma;y)
\geq \lambda_1|\xi|^2-\lambda_1^{-1}s^2\tau^2
\geq \frac{\lambda_1}{4}|\xi|^2,
\end{equation*}
that holds for  $s_0$ is sufficiently small, we can derive the estimates
\begin{equation}\label{4.82}
\begin{cases}
\begin{array}{l}
\Re \zeta_{\pm}\geq \frac{\lambda_2^2}{4}|\xi|^2,\\
R_{\pm}\geq \frac{\lambda_2}{2}|\xi|,\\
|J_{\pm}|\leq 4\lambda_2^{-3}s\tau,\\
|\partial_y \zeta_{\pm}|\leq M_4\left(1+\frac{4s_0^2\lambda_2^{-2}}{(\alpha_++L\alpha_-)^2}\right)|\xi|^2:=M_6|\xi|^2,\\
|\partial_y \sqrt{\zeta_{\pm}}|\leq \frac{M_6}{\lambda_2}|\xi|:=M_7|\xi|.
\end{array}
\end{cases}
\end{equation}

\begin{lemma}\label{lem4.7} Assume \eqref{case3}.
There exist a positive constant $C$ such that, if $s_0\le C^{-1}$ and $\tau\geq C$, then we have
\begin{equation}\label{4.83}
\Lambda|\omega_+(\xi,0)|^2+\Lambda^2\int_0^\infty|\omega_+(\xi,y)|^2dy+\int_0^\infty|\partial_y\omega_+(\xi,y)|^2dy
\leq C ||E_{+}\omega_{+}(\xi,\cdot)||^2_{L^2(\mathbb{R}_{+})}.
\end{equation}
Furthermore, if ${\rm supp}(\hat{v}_-(\xi,\cdot))\subset[-\frac{\alpha_-}{2\beta},0]$, then
\begin{equation}\label{4.84}
\Lambda^2\int_{-\infty}^0|\hat{v}_-(\xi,y)|^2dy+\int_{-\infty}^0|\partial_y\hat{v}_-(\xi,y)|^2dy
\leq C||E_{-}\hat{v}_{-}(\xi,\cdot)||^2_{L^2(\mathbb{R}_{-})}+C\Lambda|\hat{v}_-(\xi,0)|^2.
\end{equation}
\end{lemma}
\pf. 
We write
$$E_{+}\omega_{+}=I_{11}-I_{12},$$
where
\begin{eqnarray*}
I_{11}&=&\partial_y\omega_++it_{+}(\xi,y)\omega_+-iJ_+\omega_+,\\
I_{12}&=&\tau\alpha_{+}\omega_++\tau\beta y\omega_++R_+\omega_++t_{+}(\tau s\gamma,y)\omega_+,
\end{eqnarray*}
and thus
\begin{equation}\label{4.85}
\begin{aligned}
&||E_{+}\omega_{+}(\xi,\cdot)||^2_{L^2(\mathbb{R}_{+})}\\
=&\int_0^\infty|I_{11}|^2dy+\int_0^\infty[\tau\alpha_{+}+\tau\beta y+R_{+}+t_{+}(\tau s\gamma,y)]^2|\omega_+(\xi,y)|^2dy-2\Re \int_0^\infty I_{11}\bar{I}_{12}dy.
\end{aligned}
\end{equation}
We first estimate
\begin{equation}\label{4.86}
\begin{aligned}
&-2\Re\int_0^\infty I_{11}\bar{I}_{12}\\
&=\int_0^\infty[\tau\beta +\partial_yR_{+}(\xi,y)+\partial_yt_{+}(\tau s\gamma,y)]|\omega_+(\xi,y)|^2dy\\
&\quad+[\tau\alpha_{+}+R_{+}(\xi,0)+t_{+}(\tau s\gamma,0)]|\omega_+(\xi,0)|^2\\
& \geq-(M_7|\xi|-M_3s\tau)\int_0^\infty|\omega_+(\xi,y)|^2dy+\left(\tau\alpha_++\frac{\lambda_2}{\sqrt{2}}-\lambda_3^{-1}s\tau\right)|\omega_+(\xi,0)|^2
\\&\geq-C\Lambda\int_0^\infty|\omega_+(\xi,y)|^2dy+C\Lambda|\omega_+(\xi,0)|^2,
\end{aligned}
\end{equation}
provided $s_0$ is small enough. Combining \eqref{4.85} and arguing as in \eqref{4.28}, we get \eqref{4.83}.
Likewise, we obtain \eqref{4.84}.
\eproof

\begin{lemma}\label{lem4.8} Assume \eqref{case3}. There exists a positive constants $C$, depending on $\lambda_0,M_0$, such that if $s_0\le C^{-1}$, $\tau\ge C$, and $\beta\ge C$, then, for ${\rm supp}(\hat{v}_+(\xi,\cdot))\subset[0,\frac{1}{\beta}]$, we have that
\begin{equation}\label{4.89}
\frac{\Lambda^2}{\tau}\int^{\infty}_0|\hat{v}_+(\xi,y)|^2dy+\frac{1}{\tau}\int_0^{\infty}|\partial_y\hat{v}_+(\xi,y)|^2dy
\leq C\left(||F_{+}\hat{v}_{+}(\xi,\cdot)||^2_{L^2(\mathbb{R}_{+})}+\Lambda|\hat{v}_+(\xi,0)|^2\right).
\end{equation}
\end{lemma}
\pf. Expressing
$$F_{+}\hat{v}_{+}=I_{13}-I_{14},$$
where
\begin{eqnarray*}I_{13}&=&\partial_y\hat{v}_++it_{+}(\xi,y)v_++iJ_+\hat{v}_+\\
I_{14}&=&\tau\alpha_{+}\hat{v}_++\tau\beta y\hat{v}_+-R_+\hat{v}_++t_{+}(\tau s\gamma,y)\hat{v}_+,
\end{eqnarray*}
we can compute
\begin{equation}\label{4.92}
\begin{aligned}
&||F_{+}\hat{v}_{+}(\xi,\cdot)||^2_{L^2(\mathbb{R}_{+})}\\
=&\int_0^\infty|I_{13}|^2dy+\int_0^\infty p|\hat{v}_+(\xi,y)|^2dy+[\tau\alpha_{+}-R_{+}(\xi,0)+t_{+}(\tau s\gamma,0)]|\hat{v}_+(\xi,0)|^2.
\end{aligned}
\end{equation}
where $p=[-\tau\alpha_{+}-\tau\beta y+R_{+}-t_{+}(\tau s\gamma,y)]^2+(\tau\beta -\partial_yR_{+}+\partial_yt_{+}(\tau s\gamma,y))$.

We want to claim that
\begin{equation}\label{4.93}
p\geq C\frac{\Lambda^2}{\tau}.
\end{equation}
It follows from \eqref{4.82} and \eqref{case3} that for $0\leq y\leq 1/\beta$
\begin{equation}\label{4.94}
\begin{aligned}
&R_{+}-\tau\alpha_{+}-\tau\beta y-t_{+}(\tau s\gamma,y)\\
\geq& \frac{\lambda_2}{2}|\xi|-\tau(\alpha_{+}+1+\lambda_3^{-1}s_0)\geq \frac{\lambda_2}{4}|\xi|
\end{aligned}
\end{equation}
provided $|\xi|\geq C_2\tau=4\lambda_2^{-1}(\alpha_{+}+1+\lambda_3^{-1}s_0)\tau$. By \eqref{4.94}, we can easily obtain \eqref{4.93} in the case of  $|\xi|\geq C_2\tau$ with $\tau$ large.
On the other hand, when $|\xi|\leq C_2\tau$, we can estimate
\begin{equation}\label{4.95}
p\geq \tau\beta -\partial_yR_{+}+\partial_yt_{+}(\tau s\gamma,y)\geq \tau\beta-M_7C_2\tau-M_3s\tau\geq \frac{\beta}{2}\tau\geq \frac{\beta}{2}\frac{\Lambda^2}{\tau}
\end{equation}
provided $\beta$ is big enough. The estimate \eqref{4.89} is an easy consequence of \eqref{4.92} and \eqref{4.93}.
\eproof

\begin{lemma}\label{lem4.9}
Assume \eqref{case3}. There exist positive constants $C$ and $\rho_1$, depending only $\lambda_0$ and $M_0$ such that if ${\rm supp}(\hat{v}_-(\xi,\cdot))\subset[-\rho_1,0]$ then
\begin{equation}\label{4.96}
\Lambda|\omega_-(\xi,0)|^2+\Lambda^2||\omega_-(\xi,\cdot)||^2_{L^2(\mathbb{R}_{-})}
\leq C||F_{-}\omega_{-}(\xi,\cdot)||^2_{L^2(\mathbb{R}_{-})}.
\end{equation}
\end{lemma}
\pf. From \eqref{omega}, we see that
\[
{\rm supp}(\omega_-(\xi,\cdot))\subseteq{\rm supp}(\hat{v}_-(\xi,\cdot)).
\]
We first compute
\begin{equation}\label{4.97}
\begin{aligned}
&\Re\int_{-\infty}^0|\xi|(F_{-}\omega_{-})\bar{\omega}_{-}dy\\
=&\Re\int_{-\infty}^0|\xi|\partial_y\omega_-\bar{\omega}_{-}dy-\int_{-\infty}^0|\xi|[\tau\alpha_{-}+\tau\beta y+t_{-}(\tau s\gamma,y)-R_-(\xi,y)]|\omega_-|^2dy\\
=&\frac{1}{2}|\xi||\omega_-(\xi,0)|^2+\int_{-\infty}^0|\xi|[R_-(\xi;y)-\tau\alpha_{-}-\tau\beta y-t_{-}(\tau s\gamma,y)]|\omega_-|^2dy.
\end{aligned}
\end{equation}
We hope to show that
\begin{equation}\label{4.98}
C_*|\xi|\leq R_-(\xi,y)-\tau\alpha_{-}-\tau\beta y-t_{-}(\tau s\gamma,y).
\end{equation}
Assume that \eqref{4.98} is true. From \eqref{4.97} and \eqref{4.98}, it follows that
\begin{equation}\label{4.99}
\begin{aligned}
&\frac{1}{2}|\xi||\omega_-(\xi,0)|^2+\int_{-\infty}^0C_*|\xi|^2|\omega_-(\xi,y)|^2dy\\
\leq&\Re\int_{-\infty}^0|\xi|(F_{-}\omega_{-})\bar{\omega}_{-}\\
\leq& \frac{C_*}{2}\int_{-\infty}^0|\xi|^2|\omega_-(\xi,y)|^2dy+C\int_{-\infty}^0|F_{-}\omega_{-}(\xi,y)|^2dy,
\end{aligned}
\end{equation}
which implies \eqref{4.96}.

To establish \eqref{4.98}, we first note that, by simple calculations, 
\begin{equation*}
|m_-(\xi,0)-R_-(\xi,0)|\le Cs|\xi|,
\end{equation*}
which can be used to derive for $y\le 0$
\begin{equation}\label{4.100}
\begin{aligned}
&R_-(\xi,y)-\tau\alpha_{-}-\tau\beta y-t_{-}(\tau s\gamma,y)\\
\geq& m_-(\xi,0)-|R_-(\xi,0)-m_-(\xi,0)|-|R_-(\xi,y)- R_-(\xi,0)|-\tau\alpha_{-}-\lambda_3^{-1}\tau s\\
\geq & m_-(\xi,0)-\tau\alpha_{-}-C(s+|y|)|\xi|.
\end{aligned}
\end{equation}
On the other hand, by the definition of $L$, \eqref{3.6} and \eqref{case3}, we can estimate
\begin{equation}\label{4.101}
m_-(\xi,0)-\tau\alpha_{-}
\geq \frac{m_+(\xi,0)}{L}[1-\frac{L\alpha_-}{(1-\kappa)\alpha_+}]
\geq m_+(\xi,0)\frac{\kappa}{L(1-\kappa)}
\geq \frac{\lambda_2\kappa}{(1-\kappa)L} |\xi|.
\end{equation}
Combining \eqref{4.100} and \eqref{4.101} yields \eqref{4.98} provided $s$ and $|y|$ are small.
\eproof.

\begin{lemma}\label{lem4.10}
Assume \eqref{case3}. There exists $C$, depending only on $\lambda_0$ and $M_0$, such that if $s_0\le C^{-1}$, $\tau\ge C$, $\beta\ge C$, then for ${\rm supp}(\hat{v}_+(\xi,\cdot))\subset[0,\frac{1}{\beta}]$ we have
\begin{equation}\label{4.102}
\begin{aligned}
&\Lambda|V_+(\xi)+a_{nn}^+(0)\sqrt{\zeta_{+}(\xi,0)}v_+(\xi,0)|^2+\Lambda^2||F_{+}\hat{v}_{+}(\xi,\cdot)||^2_{L^2(\mathbb{R}_{+})}\\
\leq&C\left(||P_{+}\hat{v}_{+}(\xi,\cdot)||^2_{L^2(\mathbb{R}_{+})}+\Lambda^2||\hat{v}_{+}(\xi,\cdot)||^2_{L^2(\mathbb{R}_{+})}\right).
\end{aligned}
\end{equation}
Furthermore, if ${\rm supp}(\hat{v}_-(\xi,\cdot))\subset[-\rho_1,0]$, for $\rho_1$ as in Lemma \ref{lem4.9}, then
\begin{equation}\label{4.103}
\begin{aligned}
&\Lambda|V_-(\xi)-a_{nn}^-(0)\sqrt{\zeta_{-}(\xi,0)}\hat{v}_-(\xi,0)|^2+\Lambda^2||E_{-}\hat{v}_{-}(\xi,\cdot)||^2_{L^2(\mathbb{R}_{-})}\\
\leq&C\left(||P_{-}\hat{v}_{-}(\xi,\cdot)||^2_{L^2(\mathbb{R}_{-})}+\Lambda^2||\hat{v}_{-}(\xi,\cdot)||^2_{L^2(\mathbb{R}_{-})}\right).
\end{aligned}
\end{equation}
\end{lemma}
\pf.
Inequality \eqref{4.102} follows from \eqref{4.83} and \eqref{4.51}. Similarly, \eqref{4.103} follows from \eqref{4.96} and \eqref{4.52}.
\eproof

\begin{lemma}\label{lem4.11}
There exist $C$, $\rho_2$, depending only on $\lambda_0$ and $M_0$, such that if $s_0\le C^{-1}$, $\tau\ge C$, $\beta\ge C$ then for ${\rm supp}(\hat{v}_\pm(\xi,\cdot))\subset[-\rho_2,\rho_2]$  we have that \begin{equation}\label{4.105}
\begin{aligned}
&\Lambda\sum_{\pm}|V_{\pm}(\xi)|^2+\Lambda^3\sum_{\pm}|\hat{v}_{\pm}(\xi,0)|^2+\frac{\Lambda^4}{\tau}\sum_{\pm}||\hat{v}_{\pm}(\xi,\cdot)||^2_{L^2(\mathbb{R}_{\pm})}
+\frac{\Lambda^2}{\tau}\sum_{\pm}||\partial_y\hat{v}_{\pm}(\xi,\cdot)||^2_{L^2(\mathbb{R}_{\pm})}\\
\leq&C\left(\sum_{\pm}||P_{\pm}\hat{v}_{\pm}(\xi,\cdot)||^2_{L^2(\mathbb{R}_{\pm})}+\Lambda|\eta_1(\xi)|^2+\Lambda^3|\eta_0(\xi)|^2\right).
\end{aligned}
\end{equation}
\end{lemma}
\pf.
In view of \eqref{4.82},
\begin{equation*}
\begin{aligned}
&\Lambda|a_{nn}^+(0)\sqrt{\zeta_{+}(\xi,0)}\hat{v}_+(\xi,0)+a_{nn}^-(0)\sqrt{\zeta_{-}(\xi,0)}\hat{v}_+(\xi,0)|^2\\
\geq &\Lambda|a_{nn}^+(0)R_+(\xi,0)\hat{v}_+(\xi,0)+a_{nn}^-(0)R_-(\xi,0)\hat{v}_+(\xi,0)|^2\\
\geq &\frac{1}{C}\Lambda^3|\hat{v}_+(\xi,0)|^2,
\end{aligned}
\end{equation*}
hence, by \eqref{2.18} and \eqref{2.20.1}, we have
\begin{equation}\label{4.107}
\begin{aligned}
&\frac{1}{C}\Lambda^3|\hat{v}_+(\xi,0)|^2\\
\leq& \Lambda|a_{nn}^+(0)\sqrt{\zeta_{+}}\hat{v}_+(\xi,0)+a_{nn}^-(0)\sqrt{\zeta_{-}}(\hat{v}_-(\xi,0)+\eta_0)|^2\\
= &\Lambda|V_++a_{nn}^+(0)\sqrt{\zeta_{+}}\hat{v}_+(\xi,0)+a_{nn}^-(0)\sqrt{\zeta_{-}}\hat{v}_-(\xi,0)-V_--\eta_1-a_{nn}^-\sqrt{\zeta_{-}}\eta_0|^2\\
\leq& 4\left(\Lambda|V_++a_{nn}^+(0)\sqrt{\zeta_{+}}\hat{v}_+(\xi,0)|^2+\Lambda|V_--a_{nn}^-(0)\sqrt{\zeta_{-}}\hat{v}_-(\xi,0)|^2+\Lambda|\eta_1|^2+\Lambda^3|\eta_0|^2\right).
\end{aligned}
\end{equation}
From  \eqref{4.102}, \eqref{4.103} and \eqref{2.18}, we can estimate
\begin{equation}\label{stella}
\begin{aligned}
&\Lambda^3\sum_{\pm}|\hat{v}_{\pm}(\xi,0)|^2\\
\leq &C\left(\sum_{\pm}||P_{\pm}\hat{v}_{\pm}(\xi,\cdot)||^2_{L^2(\mathbb{R}_{\pm})}+\Lambda^2||\hat{v}_{+}(\xi,\cdot)||^2_{L^2(\mathbb{R}_{+})}+\Lambda|\eta_1(\xi)|^2+\Lambda^3|\eta_0(\xi)|^2\right).
\end{aligned}
\end{equation}
Again from \eqref{4.102} and \eqref{4.103}, we have the following estimate
\begin{equation}\label{4.108}
\begin{aligned}
&\Lambda|V_+|^2\\
\leq &2\Lambda|V_++a_{nn}^+(0)\sqrt{\zeta_{+}}\hat{v}_+(\xi,0)|^2+2\Lambda|a_{nn}^+(0)\sqrt{\zeta_{+}}\hat{v}_+(\xi,0)|^2\\
\leq &2\Lambda|V_++a_{nn}^+(0)\sqrt{\zeta_{+}}\hat{v}_+(\xi,0)|^2+ C\Lambda^3|\hat{v}_+(\xi,0)|^2\\
\leq& 2\Lambda|V_++a_{nn}^+(0)\sqrt{\zeta_{+}}\hat{v}_+(\xi,0)|^2+C\left(\Lambda|V_--a_{nn}^-(0)\sqrt{\zeta_{-}}\hat{v}_-(\xi,0)|^2+\Lambda|\eta_1(\xi)|^2+\Lambda^3|\eta_0(\xi)|^2\right)\\
\leq &C\left(\sum_{\pm}||P_{\pm}\hat{v}_{\pm}(\xi,\cdot)||^2_{L^2(\mathbb{R}_{\pm})}+\Lambda^2||\hat{v}_{+}(\xi,\cdot)||^2_{L^2(\mathbb{R}_{+})}+\Lambda|\eta_1(\xi)|^2+\Lambda^3|\eta_0(\xi)|^2\right).
\end{aligned}
\end{equation}
We then obtain from \eqref{2.20.1} that
\begin{equation}\label{doppiastella}
	\Lambda\sum_{\pm}|V_{\pm}(\xi)|^2\leq C\left(\sum_{\pm}||P_{\pm}\hat{v}_{\pm}(\xi,\cdot)||^2_{L^2(\mathbb{R}_{\pm})}+\Lambda^2||\hat{v}_{+}(\xi,\cdot)||^2_{L^2(\mathbb{R}_{+})}+\Lambda|\eta_1(\xi)|^2+\Lambda^3|\eta_0(\xi)|^2\right).
\end{equation}

Combining \eqref{4.84}, \eqref{4.89}, \eqref{4.102} \eqref{4.103} and \eqref{stella}, we deduce that
\begin{equation}\label{4.110}
\begin{aligned}
&\frac{\Lambda^4}{\tau}\sum_{\pm}||\hat{v}_{\pm}(\xi,\cdot)||^2_{L^2(\mathbb{R}_{\pm})}+\frac{\Lambda^2}{\tau}\sum_{\pm}||\partial_y\hat{v}_{\pm}(\xi,\cdot)||^2_{L^2(\mathbb{R}_{\pm})}\\
\leq& C\left(\Lambda^2||F_{+}\hat{v}_{+}(\xi,\cdot)||^2_{L^2(\mathbb{R}_{+})}+\Lambda^2||E_{-}\hat{v}_{-}(\xi,\cdot)||^2_{L^2(\mathbb{R}_{-})}+\Lambda^3\sum_{\pm}|\hat{v}_{\pm}(\xi,0)|^2\right)\\
\leq&C\left(\sum_{\pm}||P_{\pm}\hat{v}_{\pm}(\xi,\cdot)||^2_{L^2(\mathbb{R}_{\pm})}+\Lambda^2||\hat{v}_{+}(\xi,\cdot)||^2_{L^2(\mathbb{R}_{+})}+\Lambda|\eta_1(\xi)|^2+\Lambda^3|\eta_0(\xi)|^2\right)
\end{aligned}
\end{equation}
Finally, putting together \eqref{stella}, \eqref{doppiastella} and \eqref{4.110} yields
\begin{equation*}
\begin{aligned}
&\Lambda\sum_{\pm}|V_{\pm}|^2+\Lambda^3\sum_{\pm}|\hat{v}_{\pm}(\xi,0)|^2
+\frac{\Lambda^2}{\tau}\sum_{\pm}||\partial_y\hat{v}_{\pm}(\xi,\cdot)||^2_{L^2(\mathbb{R}_{\pm})}
+\frac{\Lambda^4}{\tau}\sum_{\pm}||\hat{v}_{\pm}(\xi,\cdot)||^2_{L^2(\mathbb{R}_{\pm})}\\
\leq&C\left( \sum_{\pm}||P_{\pm}\hat{v}_{\pm}(\xi,\cdot)||^2_{L^2(\mathbb{R}_{\pm})}+\Lambda|\eta_1|^2
+\Lambda^3|\eta_0|^2+\Lambda^2\sum_{\pm}||\hat{v}_{\pm}(\xi,\cdot)||^2_{L^2(\mathbb{R}_{\pm})}\right)
\end{aligned}
\end{equation*}
that gives \eqref{4.105} if we take $\tau$ large enough to absorb the term $C\Lambda^2\sum_{\pm}||\hat{v}_{\pm}(\xi,\cdot)||^2_{L^2(\mathbb{R}_{\pm})}$.
\eproof

Now are ready to finish the proof of Theorem~\ref{bigthm}. Combining all cases \eqref{4.377}, \eqref{4.80}, \eqref{4.105}, we conclude that
\begin{equation}\label{4.112}
\begin{aligned}
&\Lambda\sum_{\pm}|V_{\pm}|^2+\Lambda^3\sum_{\pm}|\hat{v}_{\pm}(\xi,0)|^2
+\frac{\Lambda^2}{\tau}\sum_{\pm}||\partial_y\hat{v}_{\pm}(\xi,\cdot)||^2_{L^2(\mathbb{R}_{\pm})}
+\frac{\Lambda^4}{\tau}\sum_{\pm}||\hat{v}_{\pm}(\xi,\cdot)||^2_{L^2(\mathbb{R}_{\pm})}\\
\leq&C\left(\sum_{\pm}||P_{\pm}\hat{v}_{\pm}(\xi,\cdot)||^2_{L^2(\mathbb{R}_{\pm})}+\Lambda|\eta_1|^2+\Lambda^3|\eta_0|^2\right).
\end{aligned}
\end{equation}
Recall that
\begin{equation*}
\begin{aligned}
P_{\pm}\hat{v}_{\pm}
=&\big(\partial_y-\tau\varphi_{\pm}'+it_{\pm}(\xi+i\tau s\gamma,y)\big)a_{nn}^{\pm}(y)\big(\partial_y-\tau\varphi_{\pm}'+it_{\pm}(\xi+i\tau s\gamma,y)\big)\hat{v}_{\pm}\\
&-a_{nn}^{\pm}(y)\zeta_{\pm}(\xi,y)\hat{v}_{\pm},
\end{aligned}
\end{equation*}
which implies
\begin{equation*}
|\partial^2_y\hat{v}_{\pm}|\leq C\left( |P_{\pm}\hat{v}_{\pm}|+\Lambda|\partial_y\hat{v}_{\pm}|+\Lambda^2|\hat{v}_{\pm}|\right),
\end{equation*}
where $C$ depends only on $\lambda_0$ and $M_0$. Therefore, we can derive
\begin{equation}\label{4.114}
\begin{aligned}
&\frac{1}{\tau}\sum_{\pm}||\partial^2_y\hat{v}_{\pm}(\xi,\cdot)||^2_{L^2(\mathbb{R}_{\pm})}\\
\leq&C\left( \sum_{\pm}||P_{\pm}\hat{v}_{\pm}(\xi,\cdot)||^2_{L^2(\mathbb{R}_{\pm})}+\frac{\Lambda^4}{\tau}\sum_{\pm}||\hat{v}_{\pm}(\xi,\cdot)||^2_{L^2(\mathbb{R}_{\pm})}
+\frac{\Lambda^2}{\tau}\sum_{\pm}||\partial_y\hat{v}_{\pm}(\xi,\cdot)||^2_{L^2(\mathbb{R}_{\pm})}\right).
\end{aligned}
\end{equation}
The estimate \eqref{4.4} follows directly from \eqref{4.112} and \eqref{4.114}.\eproof

\section{Step 2 - The Carleman estimate for general coefficients}

Having at disposal the Carleman estimate when $A_\pm=A_\pm(y)$, we want to derive it for $A_\pm(x,y)$.
The main idea is to "approximate" $A_\pm(x,y)$ with coefficients depending on $y$ only.
For this purpose we will make use of a special kind of partition of unity introduced in the next section.

\subsection{Partition of unity and auxiliary results}\label{preliminary}
In this section we collect some results on a partition of unity that will be crucial in our proof. In particular we will carefully describe how this partition of unity behaves with respect to the function spaces that we use.

Let $\vartheta_0\in C^\infty_0(\mathbb{R})$ and $0\leq\vartheta\leq 1$ such that
\begin{equation*}
\vartheta_0(t)=
\begin{cases}
\begin{array}{l}
1,\quad |t|\leq 1,\\
0, \quad |t|\geq 3/2.
\end{array}
\end{cases}
\end{equation*}
Also let  $\vartheta(x)=\vartheta_0(x_1)\cdots\vartheta_0(x_{n-1})$, then we have
\begin{equation*}
\vartheta(x)=
\begin{cases}
1,\quad x\in Q_1(0),\\
0, \quad x\in\mathbb{R}^{n-1}\setminus Q_{3/2}(0).
\end{cases}
\end{equation*}
Given $\mu\geq 1$ and $g\in \mathbb{Z}^{n-1}$, we define
$$x_g=g/\mu$$
and
$$\vartheta_{g,\mu}(x)=\vartheta(\mu(x-x_g)).$$
It is not hard to see that
\begin{equation*}
{\rm supp}\,\vartheta_{g,\mu}\subset Q_{3/2\mu}(x_g)\subset Q_{2/\mu}(x_g)
\end{equation*}
and
\begin{equation}\label{6.4}
|D^k\vartheta_{g,\mu}|\leq C_1\mu^k(\chi_{ Q_{3/2\mu}(x_g)}-\chi_{Q_{1/\mu}(x_g)}),\quad k=0,1,2,
\end{equation}
where $C_1\geq 1$ depends only on $n$.

For $g\in \mathbb{Z}^{n-1}$, let
$A_g=\{g'\in \mathbb{Z}^{n-1}\,|\,{\rm supp}\,\vartheta_{g',\mu}\cap{\rm supp}\,\vartheta_{g,\mu}\neq\emptyset\}$,
then
\begin{equation}\label{cardAg}
card(A_g)\mbox{  depends  only on } n.\end{equation}
 Thus, we can define
\begin{equation}\label{6.5}
\bar{\vartheta}_{\mu}(x):=\sum_{g\in \mathbb{Z}^N}\vartheta_{g,\mu}\geq 1,\quad x\in \mathbb{R}^{n-1}.
\end{equation}
It is clear that \eqref{6.4} implies
\begin{equation}\label{6.6}
|D^k\bar{\vartheta}_{\mu}|\leq C_2\mu^k,
\end{equation}
where $C_2\geq 1$ depends on $n$. If we set
$$\eta_{g,\mu}(x)=\vartheta_{g,\mu}(x)/\bar{\vartheta}_{\mu}(x),\quad x\in \mathbb{R}^{n-1},$$
then we have that
\begin{equation}\label{6.7}
\begin{cases}
\sum_{g\in \mathbb{Z}^{n-1}}\eta_{g,\mu}= 1,\quad x\in \mathbb{R}^{n-1},\\
{\rm supp}\,\eta_{g,\mu}\subset Q_{3/2\mu}(x_g)\subset Q_{2/\mu}(x_g),\\
|D^k\eta_{g,\mu}|\leq C_3\mu^k\chi_{Q_{3/2\mu}(x_g)},\quad k=0,1,2,
\end{cases}
\end{equation}
where $C_3\geq 1$ depends on $n$. 
%

\bigskip

In Section~\ref{pre} we have recalled the definition of $H^{1/2}(\mathbb{R}^{n-1})$ and its seminorm $[\cdot]_{1/2,\mathbb{R}^{n-1}}$, in what follows
we will also need the seminorm
\begin{equation}\label{semQ}
[f]_{1/2,Q_{r}}=\left[\int_{Q_{r}}\int_{Q_{r}}\frac{|f(x)-f(y)|^2}{|x-y|^n}dxdy\right]^{1/2},
\end{equation}
where $Q_{r}=Q_{r}(0)=\{x\in\mathbb{R}^{n-1}:|x_j|\leq r,\,j=1,2,\cdots ,n-1\}$.
\begin{lemma}\label{lem5.2}
Let $f\in C^\infty(\mathbb{R}^{n-1})$ and ${\rm supp} f\subset Q_{3r/4}$ for some $r\leq 1$.  There exists a positive constant $C$, depending only on $n$, such that
\begin{equation}\label{5.6}
[f]^2_{1/2,Q_{r}}+\frac{C^{-1}}{r}\int_{Q_{r}}|f(x)|^2dx\leq||f||^2_{H^{1/2}(\mathbb{R}^{n-1})}
\leq [f]^2_{1/2,Q_{r}}+\frac{C}{r}\int_{Q_{r}}|f(x)|^2dx.
\end{equation}
\end{lemma}
\pf. It follows easily from \eqref{semR} and \eqref{semQ}, that
\begin{equation}\label{5.7}
[f]^2_{1/2,\mathbb{R}^{n-1}}
=I+[f]^2_{1/2,Q_{r}},
\end{equation}
where
\begin{equation*}
I=2\int_{\mathbb{R}^{n-1}\setminus Q_{r}}\int_{Q_{r}}\frac{|f(x)-f(y)|^2}{|x-y|^n}dydx
=2\int_{\mathbb{R}^{n-1}\setminus Q_{r}}\int_{Q_{3r/4}}\frac{|f(y)|^2}{|x-y|^n}dydx.
\end{equation*}
Note that there is a positive constant $C_n<1$, depending only on $n$, such that, for $x\in \mathbb{R}^{n-1}\setminus Q_{r}$ and $y\in Q_{3r/4}$, we have
\[C_n^{-1}|x|\leq|x-y|\leq C_n|x|,\]
hence, by using Fubini theorem, there is a constant $C$ depending only on $n$, such that
\[\frac{C^{-1}}{r}\int_{Q_{r}}|f(y)|^2dy\leq I \leq\frac{C}{r}\int_{Q_{r}}|f(y)|^2dy,\]
that, together with \eqref{5.7}, gives \eqref{5.6}.
%
%
\eproof

\begin{pr}\label{pr6.1}
Let  $\{\varsigma_g\}_{g\in \mathbb{Z}^{n-1}}$ be a family of smooth functions such that ${\rm supp}\,\varsigma_{g}\subset Q_{\frac{3}{2\mu}}(x_g)$, then
\begin{equation}\label{6.8}
[\sum_{g\in \mathbb{Z}^{n-1}}\varsigma_{g}]^2_{1/2,\mathbb{R}^{n-1}}\leq C\big(\sum_{g\in \mathbb{Z}^{n-1}}[\varsigma_{g}]^2_{1/2,Q_{\frac{2}{\mu}}(x_g)}
+\sum_{g\in \mathbb{Z}^{n-1}}\mu\int_{Q_{\frac{2}{\mu}}(x_g)}|\varsigma_{g}|^2\big),
\end{equation}
where $C$ depends only on $n$.
\end{pr}
\pf. Let $x'=\mu x$ and $y'=\mu y$, then
\begin{equation*}
\begin{aligned}
&[\sum_{g}\varsigma_{g}]^2_{1/2,\mathbb{R}^{n-1}}=\int_{\mathbb{R}^{n-1}}\int_{\mathbb{R}^{n-1}}\frac{|\sum_{g}\varsigma_{g}(x)-\sum_{g}\varsigma_{g}(y)|^2}{|x-y|^n}dydx\\
=&\mu^{2-n}\int_{\mathbb{R}^{n-1}}\int_{\mathbb{R}^{n-1}}\frac{|\sum_{g}\varsigma_{g}(x'/\mu)-\sum_{g}\varsigma_{g}(y'/\mu)|^2}{|x'-y'|^n}dy'dx'.
\end{aligned}
\end{equation*}
It is enough to consider $n=2$ and denote $\varsigma_g(x'/\mu)=\varsigma_j(x')$. Note that ${\rm supp}\,\varsigma_{j}\subset Q_{\frac{3}{2}}(j)=\{x\in \mathbb{R}| |x-j|\leq 3/2\}$. We write
\begin{equation}\label{6.10}
\int_{\mathbb{R}}\int_{\mathbb{R}}\frac{|\sum_{j\in \mathbb{Z}}\varsigma_{j}(x)-\sum_{j\in \mathbb{Z}}\varsigma_{j}(y)|^2}{|x-y|^2}\,dxdy=I_1+I_2,\end{equation}
where
\[\begin{aligned}
I_1:=&\int_{\mathbb{R}^2}\frac{|\sum_{j\in \mathbb{Z}}\varsigma_{j}(x)-\sum_{j\in \mathbb{Z}}\varsigma_{j}(y)|^2}{|x-y|^2}\chi_{\{|x-y|<1\}}\,dxdy\\
I_2:=&\int_{\mathbb{R}^2}\frac{|\sum_{j\in \mathbb{Z}}\varsigma_{j}(x)-\sum_{j\in \mathbb{Z}}\varsigma_{j}(y)|^2}{|x-y|^2}\chi_{\{|x-y|\geq1\}}\,dxdy.
\end{aligned}
\]
Let us first estimate $I_2$. It is not hard to check that
\begin{equation*}
\begin{aligned}
I_2\leq&\int_{\mathbb{R}^2}\frac{2|\sum_{j\in \mathbb{Z}}\varsigma_{j}(x)|^2+2|\sum_{j\in \mathbb{Z}}\varsigma_{j}(y)|^2}{|x-y|^2}\chi_{\{|x-y|\geq1\}}dxdy\\
=&4\int_{\mathbb{R}^2}\frac{|\sum_{j\in \mathbb{Z}}\varsigma_{j}(x)|^2}{|x-y|^2}\chi_{\{|x-y|\geq1\}}dxdy
\leq8\int_{\mathbb{R}}|\sum_{j\in \mathbb{Z}}\varsigma_{j}(x)|^2dx.
\end{aligned}
\end{equation*}
Since the cardinality of $A_g$ only depends on $n$, we have for $n=2$ that
$$|\sum_{j\in \mathbb{Z}}\varsigma_{j}(x)|^2\leq 7\sum_{j\in \mathbb{Z}}|\varsigma_{j}(x)|^2$$ and hence
\begin{equation}\label{6.11}
I_2\leq56\sum_{j\in \mathbb{Z}}\int_{\mathbb{R}}|\varsigma_{j}(x)|^2dx=56\sum_{j\in \mathbb{Z}}\int_{Q_2(j)}|\varsigma_{j}(x)|^2dx.
\end{equation}
Next, for $I_1$, we can see that
\begin{equation*}
\begin{aligned}
I_1=&\int_{\mathbb{R}^2}\frac{|\sum_{j\in \mathbb{Z}}\varsigma_{j}(x)-\sum_{j\in \mathbb{Z}}\varsigma_{j}(y)|^2}{|x-y|^2}\chi_{\{|x-y|<1\}}dxdy\\
\leq&\sum_{i\in \mathbb{Z}}\int_{Q_{2}(i)}\int_{\mathbb{R}}\frac{|\sum_{j\in \mathbb{Z}}\varsigma_{j}(x)-\sum_{j\in \mathbb{Z}}\varsigma_{j}(y)|^2}{|x-y|^2}\chi_{\{|x-y|<1\}}dydx.
\end{aligned}
\end{equation*}
We note for each $x\in Q_{2}(i)$ that
$$\sum_{j\in \mathbb{Z}}\varsigma_{j}(x)=\sum_{|j-i|\leq 3}\varsigma_{j}(x)$$
and
$${\rm dist}(Q_{2}(i),Q_{2}(l))\geq 1,\quad |l-i|\geq 5.$$
Therefore, we have
\begin{equation*}
\begin{aligned}
I_1\leq&\sum_{i\in \mathbb{Z}}\int_{Q_{2}(i)}\sum_{|l-i|\leq 4}\int_{Q_{2}(l)}\frac{|\sum_{|j-i|\leq 3}\varsigma_{j}(x)-\sum_{k\in \mathbb{Z}}\varsigma_{k}(y)|^2}{|x-y|^2}dydx\\
=&\sum_{i\in \mathbb{Z}}\sum_{|l-i|\leq 4}\int_{Q_{2}(i)}\int_{Q_{2}(l)}\frac{|\sum_{|j-i|\leq 3}\varsigma_{j}(x)-\sum_{|k-l|\leq 3}\varsigma_{k}(y)|^2}{|x-y|^2}dydx.
\end{aligned}
\end{equation*}
For $|l-i|\leq 4$,  we note for $y\in Q_{2}(l)$ and $x\in Q_{2}(i)$ that
$$\sum_{|k-l|\leq 3}\varsigma_{k}(y)=\sum_{|k-i|\leq 7}\varsigma_{k}(y)$$
and
$$\sum_{|j-i|\leq 3}\varsigma_{j}(x)=\sum_{|j-i|\leq 7}\varsigma_{j}(x).$$
Thus, we can derive
\begin{equation*}
\begin{aligned}
I_1\leq&\sum_{i\in \mathbb{Z}}\int_{Q_{2}(i)}\sum_{|l-i|\leq 4}\int_{Q_{2}(l)}\frac{|\sum_{|k-i|\leq 7}(\varsigma_{k}(x)-\varsigma_{k}(y))|^2}{|x-y|^2}dydx\\
\leq&7\sum_{i\in \mathbb{Z}}\sum_{|l-i|\leq 4}\sum_{|k-i|\leq 7}\int_{Q_{2}(i)}\int_{Q_{2}(l)}\frac{|(\varsigma_{k}(x)-\varsigma_{k}(y))|^2}{|x-y|^2}dydx.
\end{aligned}
\end{equation*}
Since $Q_{2}(l)\subset Q_{6}(i)$ when $|l-i|\leq 4$, we obtain that
\begin{equation}\label{6.15}
\begin{aligned}
I_1\leq&7\sum_{i\in \mathbb{Z}}\sum_{|l-i|\leq 4}\sum_{|k-i|\leq 7}\int_{Q_{2}(i)}\int_{Q_{6}(i)}\frac{|(\varsigma_{k}(x)-\varsigma_{k}(y))|^2}{|x-y|^2}dydx\\
\leq& 10^3\sum_{k\in \mathbb{Z}}\int_{\mathbb{R}}\int_{\mathbb{R}}\frac{|(\varsigma_{k}(x)-\varsigma_{k}(y))|^2}{|x-y|^2}dydx\\
\leq& 10^3\sum_{k\in \mathbb{Z}}\{[\varsigma_{k}]^2_{1/2,Q_{2}(k)}+\int_{Q_{2}(k)}|\varsigma_{k}(x)|^2dx\},
\end{aligned}
\end{equation}
where we used \eqref{5.6} in the last inequality. Combining \eqref{6.11} and \eqref{6.15}, the proof is complete.
\eproof


\begin{pr}\label{pr6.4}
Let  $F\in C^\infty(\mathbb{R}^{n-1})\cap H^{1/2}(\mathbb{R}^{n-1})$ with ${\rm supp}\, F\subset Q_{3/2\mu}(x_g)$, and let $a$ be a function satisfying
\begin{equation}\label{6.25}
|a(z)|\leq E_a,\quad
|a(x)-a(x')|\leq K_a|x-x'|,
\end{equation}
for $z,x,x'\in {\rm supp}\,\eta_{g,\mu}$ and $E_a$, $K_a$  positive constants. Then, there is a constant $C$ depending only on $n$ such that,
\begin{equation}\label{6.26}
[a F]^2_{1/2,Q_{\frac{2}{\mu}}(x_g)}\leq C\left( E_a^2[F]^2_{1/2,Q_{\frac{2}{\mu}}(x_g)}+K_a^2\mu^{-1}\int_{Q_{\frac{2}{\mu}}(x_g)}|F(y)|^2dy\right).
\end{equation}
\end{pr}
\pf. In view of \eqref{semQ}, we have
\begin{equation*}
\begin{aligned}
&[a F]^2_{1/2,Q_{\frac{2}{\mu}}(x_g)}
=\int_{Q_{\frac{2}{\mu}}(x_g)}\int_{Q_{\frac{2}{\mu}}(x_g)}\frac{|a(x)F(x)-a(y)F(y)|^2}{|x-y|^n}\,dxdy\\
&\leq2\int_{Q_{\frac{2}{\mu}}(x_g)}\int_{Q_{\frac{2}{\mu}}(x_g)}\left(\frac{|a(x)|^2\,|F(x) -F(y)|^2}{|x-y|^{n}}+\frac{|F(y)|^2\,| a(x)-a(y)|^2}{|x-y|^{n}}\right)dxdy,\\
&\leq C\left(E_a^2[F]^2_{1/2,Q_{\frac{2}{\mu}}(x_g)}+K_a^2\mu^{-1}\int_{Q_{\frac{2}{\mu}}(x_g)}|F(y)|^2dy\right).
\end{aligned}
\end{equation*}
\eproof

\begin{pr}\label{pr6.2}
Let  $f\in C^\infty(\mathbb{R}^{n-1})\cap H^{1/2}(\mathbb{R}^{n-1})$. Then
\begin{equation}\label{6.16}
\sum_{g\in \mathbb{Z}^{n-1}}[f \eta_{g,\mu}]^2_{1/2,Q_{\frac{2}{\mu}}(x_g)}\leq C\left( [f]^2_{1/2,\mathbb{R}^{n-1}}+\mu\int_{\mathbb{R}^{n-1}}|f(y)|^2dy\right).
\end{equation}
\end{pr}
\pf. It follows from \eqref{semQ} that
\begin{equation}\label{6.17}
[f \eta_{g,\mu}]^2_{1/2,Q_{\frac{2}{\mu}}(x_g)}
\leq I+2\int_{Q_{\frac{2}{\mu}}(x_g)}\int_{Q_{\frac{2}{\mu}}(x_g)}\frac{|\eta_{g,\mu}(x)|^2\,|f(x) -f (y)|^2}{|x-y|^{n}}\,dxdy,
\end{equation}
where $$I=2\int_{Q_{\frac{2}{\mu}}(x_g)}\int_{Q_{\frac{2}{\mu}}(x_g)}|f(y)|^2\,| \eta_{g,\mu}(x)-\eta_{g,\mu}(y)|^2\,|x-y|^{-n}dxdy.$$
Using \eqref{6.7}, we can estimate
\begin{equation}\label{6.18}
I\leq 2C_3^2\mu^{2}\int_{Q_{\frac{2}{\mu}}(x_g)}\int_{Q_{\frac{2}{\mu}}(x_g)}
\frac{|f(y)|^2}{|x-y|^{n}}\,dxdy\\
\leq C\mu\int_{Q_{\frac{2}{\mu}}(x_g)}|f(y)|^2dy.
\end{equation}
Using the fact that for any $g\in\mathbb{Z}^{n-1}$, the cardinality of
$\{g'\in\mathbb{Z}^{n-1}| Q_{2/\mu}(x_g)\cap Q_{2/\mu}(x_{g'})\ne\emptyset\}$ is finite and only depends on $n$ and adding up with respect to $g\in \mathbb{Z}^{n-1}$, we get \eqref{6.16}.
\eproof

\begin{pr}\label{pr6.3}
Let  $f\in C^\infty(\mathbb{R}^{n-1})\cap H^{1/2}(\mathbb{R}^{n-1})$. Then
\begin{equation}\label{6.20}
\sum_{g\in \mathbb{Z}^{n-1}}[f\, \nabla_x\eta_{g,\mu}]^2_{1/2,Q_{\frac{2}{\mu}}(x_g)}\leq C\left(\mu^2[f]^2_{1/2,\mathbb{R}^{n-1}}+\mu^3\int_{\mathbb{R}^{n-1}}|f(y)|^2dy\right).
\end{equation}
\end{pr}
We omit the proof that proceeds in the same way as that of Proposition~\ref{pr6.2}.

\subsection{Estimate of the left hand side of the Carleman estimate, I}

We are ready to derive the Carleman estimate for general coefficients. In order to make clear the procedure that we follow let us introduce and recall some notations and definitions. Let $0<\delta\leq 1$ and define
\begin{equation}\label{7.7}
A^\delta_{\pm}(x,y):=A_{\pm}(\delta x,\delta y),
\end{equation}
\begin{equation}\label{7.8}
\mathcal{L}_\delta(x,y,\partial)w:=\sum_{\pm}H_{\pm}{\rm div}_{x,y}(A^\delta_{\pm}(x,y)\nabla_{x,y}w_{\pm}),
\end{equation}
and the transmission conditions
\begin{equation*}
\begin{cases}
\theta_0(x)=w_+(x,0)-w_-(x,0),\\
\theta_1(x)=A_+(x,0)\nabla_{x,y}w_+(x,0)\cdot \nu-A_-(x,0)\nabla_{x,y}w_-(x,0)\cdot \nu.
\end{cases}
\end{equation*}
Next, with $x_g=g/\mu$ $g\in \mathbb{Z}^{n-1}$, we define
\begin{equation*}
\begin{cases}
A^{\delta,g}_{\pm}(y):=A^\delta_{\pm}(x_g, y)=A_{\pm}(\delta x_g,\delta y),\\
\mathcal{L}_{\delta,g}(y,\partial)w:=\sum_{\pm}H_{\pm}{\rm div}_{x,y}(A^{\delta,g}_{\pm}(y)\nabla_{x,y}w_{\pm}).
\end{cases}
\end{equation*}
It is readily seen that
\begin{equation*}
\lambda_0|z|^2\leq A^{\delta,g}_{\pm}(y)z\cdot z\leq \lambda^{-1}_0|z|^2,\, \forall y\in \mathbb{R},\,\forall\, z\in \mathbb{R}^n
\end{equation*}
and
\begin{equation*}
|A^{\delta,g}_{\pm}(y')-A^{\delta,g}_{\pm}(y)|\leq M_0\delta|y'-y|.
\end{equation*}
Concerning the weight functions, let us introduce the following notations
\[
\begin{cases}
h_\varepsilon(x):=-\varepsilon|x|^2/2,\\
H_\varepsilon(x,x_g):=\varepsilon|x-x_g|^2/2,\\
\psi_\varepsilon(x,y):=\varphi(y)+h_\varepsilon(x),\\
\psi_{\varepsilon,g}(x,y):=\varphi(y)+\nabla_xh_\varepsilon(x_g)\cdot(x-x_g)+h_\varepsilon(x_g),
\end{cases}
\]
where $\varphi(y)$ is defined in \eqref{2.1}. Moreover assume that $\alpha_{+}, \alpha_{-}, \beta$ are fixed positive numbers such that $\beta\geq\beta_0$ and $\lambda_2^{-1}<\frac{\alpha_{+}}{\alpha_{-}}$, in such a way that condition \eqref{4.3} is satisfied by the operator $\mathcal{L}_{\delta,g}(y,\partial)$ and Theorem \ref{pr22} holds true for such an operator.

Note that
\begin{equation}\label{7.9}
\psi_{\varepsilon,g}(x,y)-\psi_\varepsilon(x,y)=H_\varepsilon(x,x_g).
\end{equation}
We define
\begin{equation}\label{l.s.}
\begin{aligned}
\Xi(w):=&\sum_{\pm}\sum_{k=0}^2\tau^{3-2k}\int_{\mathbb{R}^n_{\pm}}|D^k{w}_{\pm}|^2e^{2\tau\psi_{\varepsilon}}dxdy\\
&+\sum_{\pm}\sum_{k=0}^1\tau^{3-2k}\int_{\mathbb{R}^{n-1}}|D^k{w}_{\pm}(x,0)|^2e^{2\tau\psi_{\varepsilon}(x,0)}dx\\
&+\sum_{\pm}\tau^2[e^{\tau\psi_{\varepsilon}(\cdot,0)}w_{\pm}(\cdot,0)]^2_{1/2,\mathbb{R}^{n-1}}\\
&+\sum_{\pm}[\partial_y(e^{\tau\psi_{\varepsilon,\pm}}w_{\pm})(\cdot,0)]^2_{1/2,\mathbb{R}^{n-1}}
+\sum_{\pm}[\nabla_x(e^{\tau\psi_{\varepsilon}}w_{\pm})(\cdot,0)]^2_{1/2,\mathbb{R}^{n-1}},
\end{aligned}
\end{equation}
where we note that $\Xi(w)$ corresponds to the left hand side of \eqref{8.24}.

In the present subsection we prove that if ${\rm supp}\, w\subset \mathfrak{U}:=B_{1/2}\times [-r_0,r_0]$ and if we choose
\begin{equation}\label{7.17}
\tau\geq 1/\varepsilon\quad\text{and}\quad  \mu=(\varepsilon\tau)^{1/2},
\end{equation}
then
\begin{equation}\label{7.29}
\begin{array}{l}
\Xi(w)\leq C\sum_{g\in \mathbb{Z}^{n-1}} \Xi(w\eta_{g,\mu})+CR_1,
\end{array}
\end{equation}
where
\begin{equation*}
R_1:= (\varepsilon\tau)^{1/2}\sum_{\pm}\int_{\mathbb{R}^{n-1}}e^{2\tau\psi_{\varepsilon}(x,0)}(|\partial_yw_{\pm}(x,0)|^2+|\nabla_xw_{\pm}(x,0)|^2+\tau^2|w_{\pm}(x,0)|^2)dx
\end{equation*}
and $C$ depends only on $\lambda_0,M_0$. 

In order to obtain \eqref{7.29}, we estimate from above each term in \eqref{l.s.}, by \eqref{6.7}, 
\begin{equation}\label{7.12}
w_{\pm}=\sum_{g\in \mathbb{Z}^{n-1}}w_{\pm}\eta_{g,\mu}(x).
\end{equation}
From \eqref{cardAg}, \eqref{7.9} and \eqref{7.12}, we can see that
\begin{equation}\label{7.15}
\begin{aligned}
&
\sum_{k=0}^2\tau^{3-2k}\int_{\mathbb{R}^n_{\pm}}|D^k{w}_{\pm}|^2e^{2\tau\psi_{\varepsilon}}\,dxdy\\
&\leq	C \sum_{g\in \mathbb{Z}^{n-1}}
\sum_{k=0}^2\tau^{3-2k}\int_{\mathbb{R}^n_{\pm}}|D^k(w_{\pm}\eta_{g,\mu})|^2e^{2\tau\psi_{\varepsilon,g,\pm}}dxdy
\end{aligned}
\end{equation}
and
\begin{equation}\label{7.15,5}
\begin{aligned}
&
\sum_{k=0}^1\tau^{3-2k}\int_{\mathbb{R}^{n-1}}|D^k{w}_{\pm}(x,0)|^2e^{2\tau\psi_{\varepsilon}(x,0)}\,dx\\
&\leq C \sum_{g\in \mathbb{Z}^{n-1}}
\sum_{k=0}^1\tau^{3-2k}\int_{\mathbb{R}^{n-1}}|D^k(w_{\pm}\eta_{g,\mu})(x,0)|^2e^{2\tau\psi_{\varepsilon,g}(x,0)}\,dx,\end{aligned}
\end{equation}
where $C$ depends only on $n$.

Using \eqref{6.8}, we obtain
\begin{equation}\label{7.16}
\begin{aligned}
& [\nabla_x(e^{\tau\psi_{\varepsilon}}w_{\pm})(\cdot,0)]^2_{1/2,\mathbb{R}^{n-1}}=
[\nabla_x(e^{\tau\psi_{\varepsilon}}\sum_{g}w_{\pm}\eta_{g,\mu})(\cdot,0)]^2_{1/2,\mathbb{R}^{n-1}}\\
&\leq C
\sum_{g\in \mathbb{Z}^{n-1}}\left([\nabla_x(e^{\tau\psi_{\varepsilon}}\eta_{g,\mu}w_{\pm})(\cdot,0)]^2_{1/2,Q_{\frac{2}{\mu}}(x_g)}
+
\mu\int_{Q_{\frac{2}{\mu}}(x_g)}|\nabla_x(e^{\tau\psi_{\varepsilon}}\eta_{g,\mu}w_{\pm})(x,0)|^2dx\right).
\end{aligned}
\end{equation}
Since
$$
\begin{aligned}
&\nabla_x(e^{\tau\psi_{\varepsilon}}\eta_{g,\mu}w_{\pm})(x,0)\\
=&e^{\tau\psi_{\varepsilon}(x,0)}\eta_{g,\mu}\nabla_xw_{\pm}(x,0)+e^{\tau\psi_{\varepsilon}(x,0)}w_{\pm}\nabla_x\eta_{g,\mu}(x,0)
-(\varepsilon\tau x)e^{\tau\psi_{\varepsilon}(x,0)}\eta_{g,\mu}w_{\pm}(x,0),
\end{aligned}
$$
by \eqref{6.8} and \eqref{cardAg}, we have that
\begin{equation}\label{7.18}
\begin{aligned}
&
\sum_{g\in \mathbb{Z}^{n-1}}\mu\int_{Q_{\frac{2}{\mu}}(x_g)}|\nabla_x(e^{\tau\psi_{\varepsilon}}\eta_{g,\mu}w_{\pm})(x,0)|^2dx\\
\leq& 
C\left(\mu\int_{\mathbb{R}^{n-1}}e^{2\tau\psi_{\varepsilon}(x,0)}|\nabla_xw_{\pm}(x,0)|^2dx
+
\mu^5\int_{\mathbb{R}^{n-1}}e^{2\tau\psi_{\varepsilon}(x,0)}|w_{\pm}(x,0)|^2dx\right).
\end{aligned}
\end{equation}

In order to estimate $\sum_{g\in \mathbb{Z}^{n-1}}[\nabla_x(e^{\tau\psi_{\varepsilon}}\eta_{g,\mu}w_{\pm})(\cdot,0)]^2_{1/2,Q_{2/\mu}(x_g)}$ we need to observe the following easy consequence of Lemma \ref{lem5.2}, Proposition \ref{pr6.4} and of \eqref{7.9}:

\begin{lemma}\label{lemma7.2}
If ${\rm supp}\, f\subset Q_{3/2\mu}(x_g)$, then we have that
\begin{equation}\label{7.19}
[f e^{\tau\psi_{\varepsilon}(\cdot,0)}]^2_{1/2,\mathbb{R}^{n-1}}
\leq C 
[f e^{\tau\psi_{\varepsilon,g}(\cdot,0)}]^2_{1/2,\mathbb{R}^{n-1}}+\mu\int_{Q_{2/\mu}(x_g)}|f(x) |^2e^{2\tau\psi_{\varepsilon}(x,0)}dx,
\end{equation}
and
\begin{equation}\label{7.191}
[f e^{\tau\psi_{\varepsilon,g}(\cdot,0)}]^2_{1/2,\mathbb{R}^{n-1}}
\leq C 
[f e^{\tau\psi_{\varepsilon}(\cdot,0)}]^2_{1/2,Q_{2/\mu}(x_g)}+\mu\int_{Q_{2/\mu}(x_g)}|f(x) |^2e^{2\tau\psi_{\varepsilon}(x,0)}dx,
\end{equation}
where $C$ depends only on $n$.
\end{lemma}

\noindent Similarly, we can show that
\begin{lemma}\label{lemma7.3}
\begin{equation}\label{7.23}
\begin{aligned}
&[x e^{\tau\psi_{\varepsilon}(\cdot,0)}\eta_{g,\mu}w_{\pm}]^2_{1/2,Q_{2/\mu}(x_g)}\leq C\left([e^{\tau\psi_{\varepsilon,g}(\cdot,0)}\eta_{g,\mu}w_{\pm}]^2_{1/2,Q_{2/\mu}(x_g)}\right.\\
&\left.\hspace{20mm}+\frac{1}{\mu}\int_{Q_{2/\mu}(x_g)}|\eta_{g,\mu}w_{\pm}(x,0)|^2e^{2\tau\psi_{\varepsilon}(x,0)}dx\right).
\end{aligned}
\end{equation}
\end{lemma}

It is time to estimate $\sum_{\pm}\sum_{g\in \mathbb{Z}^{n-1}}[\nabla_x(e^{\tau\psi_{\varepsilon}}\eta_{g,\mu}w_{\pm})(\cdot,0)]^2_{1/2,Q_{\frac{2}{\mu}}(x_g)}$.
Since
$$
\nabla_x(e^{\tau\psi_{\varepsilon}}\eta_{g,\mu}w_{\pm})(x,0)
=e^{\tau\psi_{\varepsilon}}\nabla_x(\eta_{g,\mu}w_{\pm})(x,0)-(\varepsilon\tau x)e^{\tau\psi_{\varepsilon}}\eta_{g,\mu}w_{\pm}(x,0),
$$
we can deduce from \eqref{7.19} and \eqref{7.23} that
\begin{equation}\label{7.25}
\begin{aligned}
&\sum_{g\in \mathbb{Z}^{n-1}}[\nabla_x(e^{\tau\psi_{\varepsilon}}\eta_{g,\mu}w_{\pm})(\cdot,0)]^2_{1/2,Q_{\frac{2}{\mu}}(x_g)}\leq \sum_{g\in \mathbb{Z}^{n-1}}[e^{\tau\psi_{\varepsilon}(\cdot,0)}\nabla_x(\eta_{g,\mu}w_{\pm})(\cdot,0)]^2_{1/2,Q_{\frac{2}{\mu}}(x_g)}\\
&+(\varepsilon\tau)^2\sum_{g\in \mathbb{Z}^{n-1}}[xe^{\tau\psi_{\varepsilon}(\cdot,0)}\eta_{g,\mu}w_{\pm})(\cdot,0)]^2_{1/2,Q_{\frac{2}{\mu}}(x_g)}\\
\leq&C\left( \sum_{g\in \mathbb{Z}^{n-1}}[e^{\tau\psi_{\varepsilon,g}(\cdot,0)}\nabla_x(\eta_{g,\mu}w_{\pm})]^2_{1/2,Q_{\frac{2}{\mu}}(x_g)}+\mu\sum_{g\in \mathbb{Z}^{n-1}}\int_{Q_{2/\mu}(x_g)}|\nabla_x(\eta_{g,\mu}w_{\pm})|^2e^{2\tau\psi_{\varepsilon}(x,0)}dx\right.\\
&\left.+(\varepsilon\tau)^2\sum_{g\in \mathbb{Z}^{n-1}}[e^{\tau\psi_{\varepsilon,g}(\cdot,0)}\eta_{g,\mu}w_{\pm})]^2_{1/2,Q_{\frac{2}{\mu}}(x_g)}+\mu^{-1}(\varepsilon\tau)^2\sum_{g\in \mathbb{Z}^{n-1}}\int_{Q_{2/\mu}(x_g)}|\eta_{g,\mu}w_{\pm}|^2e^{2\tau\psi_{\varepsilon}(x,0)}dx\right)\\
\leq &C\left(\sum_{g\in \mathbb{Z}^{n-1}}[e^{\tau\psi_{\varepsilon,g}(\cdot,0)}\nabla_x(\eta_{g,\mu}w_{\pm})]^2_{1/2,Q_{\frac{2}{\mu}}(x_g)}+(\varepsilon\tau)^2\sum_{g\in \mathbb{Z}^{n-1}}[e^{\tau\psi_{\varepsilon,g}(\cdot,0)}\eta_{g,\mu}w_{\pm})]^2_{1/2,Q_{\frac{2}{\mu}}(x_g)}\right.\\
&\left.+\sum_{\pm}\mu\int_{\mathbb{R}^{n-1}}e^{2\tau\psi_{\varepsilon}(x,0)}|\nabla_xw_{\pm}(x,0)|^2dx+\mu^3\int_{\mathbb{R}^{n-1}}e^{2\tau\psi_{\varepsilon}(x,0)}|w_{\pm}(x,0)|^2dx\right).
\end{aligned}
\end{equation}
Combining \eqref{7.16}, \eqref{7.18} and \eqref{7.25} yields
\begin{equation}\label{7.26}
\begin{aligned}
&
[\nabla_x(e^{\tau\psi_{\varepsilon}}\sum_{g}w_{\pm}\eta_{g,\mu})(\cdot,0)]^2_{1/2,\mathbb{R}^{n-1}}\leq C\left( 
\sum_{g\in \mathbb{Z}^{n-1}}[e^{\tau\psi_{\varepsilon,g}(\cdot,0)}\nabla_x(\eta_{g,\mu}w_{\pm})(\cdot,0)]^2_{1/2,Q_{\frac{2}{\mu}}(x_g)}\right.\\
&+\mu^4
\sum_{g\in \mathbb{Z}^{n-1}}[e^{\tau\psi_{\varepsilon,g}(\cdot,0)}\eta_{g,\mu}w_{\pm}(\cdot,0)]^2_{1/2,Q_{\frac{2}{\mu}}(x_g)}+
\mu\int_{\mathbb{R}^{n-1}}e^{2\tau\psi_{\varepsilon}(x,0)}|\nabla_xw_{\pm}(x,0)|^2dx\\
&+
\left.\mu^5\int_{\mathbb{R}^{n-1}}e^{2\tau\psi_{\varepsilon}(x,0)}|w_{\pm}(x,0)|^2dx\right).
\end{aligned}
\end{equation}

In a similar way, we can estimate the terms $[\partial_y(e^{\tau\psi_{\varepsilon\pm}}\sum_{g}w_{\pm}\eta_{g,\mu})(\cdot,0)]^2_{1/2,\mathbb{R}^{n-1}}$ and
$\tau^2[e^{\tau\psi_{\varepsilon}(\cdot,0)}\sum_{g}w_{\pm}\eta_{g,\mu}]^2_{1/2,\mathbb{R}^{n-1}}$
and finally get \eqref{7.29}. Notice that in deriving \eqref{7.29} we make use of $\mu^4=(\eps\tau)^2\le\tau^2$.

\subsection{Estimate of the left hand side of the Carleman estimate, II}\label{subsEII}

In this section, we will continue to estimate the upper bound of $\Xi(w)$ using \eqref{7.29}. The task now is to connect the estimate to the operator $\mathcal{L}(x,y,\partial)$ given in \eqref{7.1}. To this aim we apply Theorem \ref{pr22} to the function $w\eta_{g,\mu}$ with the weight function $\psi_{\eps,g}=\varphi(y)-\eps x_g\cdot x+\eps|x_g|^2/2$. In order to do this we note that if ${\rm supp}\, w\subset \mathfrak{U}:=B_{1/2}\times [-r_0,r_0]$ and $\mu\geq 4$, then either $|x_g|\leq 1$ or ${\rm supp}\, \eta_{g,\mu}\cap B_{1/2}=\emptyset$. Thus, in both the cases, we can apply Theorem \ref{pr22}.

By applying \eqref{5.16} and by adding up with respect to $g\in\mathbb{Z}^{n-1}$, we obtain that
\begin{equation}\label{8.1}
\sum_{g\in \mathbb{Z}^{n-1}}\Xi(w\eta_{g,\mu})\leq C \sum_{g\in \mathbb{Z}^{n-1}}(d^{(1)}_{g,\mu}+d^{(2)}_{g,\mu}+d^{(3)}_{g,\mu}),
\end{equation}
where
\begin{equation*}
\begin{aligned}
d^{(1)}_{g,\mu}=& \int_{\mathbb{R}^n}|\mathcal{L}_{\delta,g}(y,\partial)(w\eta_{g,\mu})|^2e^{2\tau\psi_{\varepsilon,g}}dxdy,\\
d^{(2)}_{g,\mu}=&\tau^{3}\int_{\mathbb{R}^{n-1}}|e^{\tau\psi_{\varepsilon,g}(x,0)}\theta_{0;g,\mu}(x)|^2dx+[\nabla_x(e^{\tau\psi_{\varepsilon,g}}\theta_{0;g,\mu})(\cdot,0)]^2_{1/2,\mathbb{R}^{n-1}},\\
d^{(3)}_{g,\mu}=&\tau\int_{\mathbb{R}^{n-1}}|e^{\tau\psi_{\varepsilon,g}(x,0)}\theta_{1;g,\mu}(x)|^2dx+[e^{\tau\psi_{\varepsilon,g}(\cdot,0)}\theta_{1;g,\mu}(\cdot)]^2_{1/2,\mathbb{R}^{n-1}},
\end{aligned}
\end{equation*}
where we set
\begin{equation}\label{theta0g}
\theta_{0;g,\mu}(x):=w_+(x,0)\eta_{g,\mu}(x)-w_-(x,0)\eta_{g,\mu}(x)=\theta_0(x)\eta_{g,\mu},
\end{equation}
\begin{equation}\label{theta1g}
\theta_{1;g,\mu}(x):=A^{\delta,g}_+(0)\nabla_{x,y}(w_+\eta_{g,\mu})\cdot \nu-A^{\delta,g}_-(0)\nabla_{x,y}(w_-\eta_{g,\mu})\cdot \nu.\end{equation}
We will estimate the three terms of \eqref{8.1} separately.

\medskip
\noindent{\bf Estimate of $\sum_{g\in \mathbb{Z}^{n-1}}d^{(1)}_{g,\mu}$}.
\medskip

By straightforward computations, we obtain that
\begin{equation*}
\begin{aligned}
&\mathcal{L}_{\delta,g}(y,\partial)(w_\pm\eta_{g,\mu})\\
\leq &|\mathcal{L}_\delta(x,y,\partial)(w_\pm\eta_{g,\mu})|+|\mathcal{L}_\delta(x,y,\partial)(w_\pm\eta_{g,\mu})-\mathcal{L}_{\delta,g}(y,\partial)(w_\pm\eta_{g,\mu})|\\
\leq &\eta_{g,\mu}|\mathcal{L}_\delta(x,y,\partial)(w_\pm)|+C\left(\delta\mu^{-1}|D^2w_\pm|\,\chi_{Q_{\frac{2}{\mu}}(x_g)}+\mu|Dw_\pm|\,\chi_{Q_{\frac{2}{\mu}}(x_g)}+\mu^2|w_\pm|\,\chi_{Q_{\frac{2}{\mu}}(x_g)}\right),
\end{aligned}
\end{equation*}
which implies
\begin{equation}\label{8.3}
\sum_{g\in \mathbb{Z}^{n-1}}d^{(1)}_{g,\mu}\leq
C\sum_{\pm}\int_{\mathbb{R}^n_{\pm}}|\mathcal{L}_\delta(x,y,\partial)(w_\pm)|^2\,e^{2\tau\psi_{\varepsilon}}dxdy+CR_2,
\end{equation}
where
\begin{equation*}
\begin{aligned}
R_2=&\delta^2\mu^{-2}\sum_{\pm}\int_{\mathbb{R}^n_{\pm}}|D^2w_{\pm}|^2\,e^{2\tau\psi_{\varepsilon,\pm}}dxdy+\mu^{2}\sum_{\pm}\int_{\mathbb{R}^n_{\pm}}|Dw_{\pm}|^2\,e^{2\tau\psi_{\varepsilon,\pm}}dxdy\\
&+\mu^{4}\sum_{\pm}\int_{\mathbb{R}^n_{\pm}}|w_{\pm}|^2\,e^{2\tau\psi_{\varepsilon,\pm}}dxdy.
\end{aligned}
\end{equation*}

\noindent{\bf Estimate of $\sum_{g\in \mathbb{Z}^{n-1}}d^{(2)}_{g,\mu}$}.
\medskip

It is obvious that
\begin{equation}\label{8.4}
\sum_{g\in \mathbb{Z}^{n-1}}\tau^{3}\int_{\mathbb{R}^{n-1}}|e^{\tau\psi_{\varepsilon,g}(x,0)}\theta_{0;g,\mu}(x)|^2dx
C\tau^{3}\int_{\mathbb{R}^{n-1}}e^{2\tau\psi_{\varepsilon}(x,0)}|\theta_0(x)|^2dx,
\end{equation}
where $C$ depends only on $n$. Next, we note that $\nabla_x(e^{\tau\psi_{\eps,g}}\theta_{0;g,\mu})=e^{\tau\psi_{\eps,g}}\nabla_x\theta_{0;g,\mu}-\tau\eps x_ge^{\tau\psi_{\eps,g}}\theta_{0;g,\mu}$. From \eqref{7.191}, \eqref{6.16}, and \eqref{theta0g}, it follows that
\begin{equation}\label{8.5}
\begin{aligned}
&\sum_{g\in \mathbb{Z}^{n-1}}[\nabla_x(e^{\tau\psi_{\varepsilon,g}}\theta_{0;g,\mu})(\cdot,0)]^2_{1/2,Q_{\frac{2}{\mu}}(x_g)}\\
\leq&C\left( \sum_{g\in \mathbb{Z}^{n-1}}[e^{\tau\psi_{\varepsilon}(\cdot,0)}\nabla_x\theta_{0;g,\mu}]^2_{1/2,Q_{\frac{2}{\mu}}(x_g)}+(\tau\eps)^2[e^{\tau\psi_{\varepsilon}(\cdot,0)}\omega_0]^2_{1/2,\mathbb{R}^{n-1}}\right.\\
&\left.+\mu\int_{\mathbb{R}^{n-1}}e^{2\tau\psi_{\varepsilon}(x,0)}|\nabla_x\theta_0|^2dx+\mu^5\int_{\mathbb{R}^{n-1}}e^{2\tau\psi_{\varepsilon}(x,0)}|\theta_0|^2dx\right).
\end{aligned}
\end{equation}
On the other hand, by \eqref{6.16}, \eqref{6.20} and \eqref{7.23}, we have that
\begin{equation}\label{8.8}
\begin{aligned}
&\sum_{g\in \mathbb{Z}^{n-1}}[e^{\tau\psi_{\varepsilon}(\cdot,0)}\nabla_x\theta_{0;g,\mu}]^2_{1/2,Q_{\frac{2}{\mu}}(x_g)}\\
\lesssim& [\nabla_x(e^{\tau\psi_{\varepsilon}(\cdot,0)}\theta_0)]^2_{1/2,\mathbb{R}^{n-1}}+\mu^4[e^{\tau\psi_{\varepsilon}(\cdot,0)}\theta_0]^2_{1/2,\mathbb{R}^{n-1}}\\
&+\mu\int_{\mathbb{R}^{n-1}}e^{2\tau\psi_{\varepsilon}(x,0)}|\nabla_x\theta_0|^2dx+\mu^5\int_{\mathbb{R}^{n-1}}e^{2\tau\psi_{\varepsilon}(x,0)}|\theta_0|^2dx.
\end{aligned}
\end{equation}
Finally, combining \eqref{8.4} and \eqref{8.8} yields
\begin{equation}\label{8.9}
\sum_{g\in \mathbb{Z}^{n-1}}d^{(2)}_{g,\mu}\\
\leq C\left( [\nabla_x(e^{\tau\psi_{\varepsilon}}\theta_0)(\cdot,0)]^2_{1/2,\mathbb{R}^{n-1}}+\tau^{3}\int_{\mathbb{R}^{n-1}}e^{2\tau\psi_{\varepsilon}(x,0)}|\theta_0|^2dx+R_3\right),
\end{equation}
where
$$
\begin{aligned}
R_3=\mu^4[e^{\tau\psi_{\varepsilon}(\cdot,0)}\theta_0]^2_{1/2,\mathbb{R}^{n-1}}+\mu\int_{\mathbb{R}^{n-1}}e^{2\tau\psi_{\varepsilon}(x,0)}|\nabla_x\theta_0|^2dx+\mu^5\int_{\mathbb{R}^{n-1}}e^{2\tau\psi_{\varepsilon}(x,0)}|\theta_0|^2dx.
\end{aligned}
$$

\medskip
\noindent{\bf Estimate of $\sum_{g\in \mathbb{Z}^{n-1}}d^{(3)}_{g,\mu}$}.
\medskip

Straightforward computations show that we can write $\theta_{1;g,\mu}$ as
\begin{equation}\label{8.10}
\theta_{1;g,\mu}=
\theta_1\eta_{g,\mu}+J^{(1)}_{g,\mu}+J^{(2)}_{g,\mu}+J^{(3)}_{g,\mu},
\end{equation}
where
$$
\begin{aligned}
J^{(1)}_{g,\mu}=&w_+A_+(\delta x,0)\nabla_{x,y}\eta_{g,\mu}\cdot \nu -w_-A_-(\delta x,0)\nabla_{x,y}\eta_{g,\mu}\cdot \nu ,\\
J^{(2)}_{g,\mu}=&\eta_{g,\mu}(A_+(\delta x_g,0)-A_+(\delta x,0))\nabla_{x,y}w_+\cdot \nu\\
&-\eta_{g,\mu}(A_-(\delta x_g,0)-A_-(\delta x,0))\nabla_{x,y}w_-\cdot \nu,\\
J^{(3)}_{g,\mu}=&w_+(A_+(\delta x_g,0)-A_+(\delta x,0))\nabla_{x,y}\eta_{g,\mu}\cdot \nu\\
&-w_-(A_-(\delta x_g,0)-A_-(\delta x,0))\nabla_{x,y}\eta_{g,\mu}\cdot \nu.
\end{aligned}
$$
It is easy to compute that
\begin{equation}\label{8.11}
\begin{aligned}
&|J^{(1)}_{g,\mu}|\leq C \mu\sum_{\pm}|w_{\pm}(x,0)|\,\chi_{Q_{\frac{2}{\mu}}(x_g)},\\
&|J^{(2)}_{g,\mu}|\leq C \delta\mu^{-1}\sum_{\pm}|\nabla_{x,y}w_{\pm}(x,0)|\eta_{g,\mu},\\
&|J^{(3)}_{g,\mu}|\leq C \delta\mu^{-1}\sum_{\pm}|\nabla_{x,y}\eta_{g,\mu}||w_{\pm}(x,0)|,
\end{aligned}
\end{equation}
where $C$ depends on $\lambda_0$, $M_0$ and $n$. Putting together \eqref{8.10} and \eqref{8.11} implies
\begin{equation}\label{8.12}
\begin{aligned}
&\sum_{g\in \mathbb{Z}^{n-1}}\tau\int_{\mathbb{R}^{n-1}}|e^{\tau\psi_{\varepsilon,g}(x,0)}\theta_{1;g,\mu}(x)|^2dx\leq C\left(\tau\int_{\mathbb{R}^{n-1}}|\omega_1|^2e^{2\tau\psi_{\varepsilon}(x,0)}dx\right.\\
&+\delta^2\varepsilon^{-1}\sum_{\pm}\int_{\mathbb{R}^{n-1}}|\nabla_{x,y}w_{\pm}(x,0)|^2e^{2\tau\psi_{\varepsilon}(x,0)}dx\\
&\left.+(\delta^{2}\tau+\tau^2\varepsilon)\sum_{\pm}\int_{\mathbb{R}^{n-1}}|w_{\pm}(x,0)|^2e^{2\tau\psi_{\varepsilon}(x,0)}dx\right).
\end{aligned}
\end{equation}

We turn to the second term of $d^{(3)}_{g,\mu}$. We first derive from \eqref{7.191} and \eqref{6.16} that
\begin{equation}\label{8.13}
\begin{aligned}
&\sum_{g\in \mathbb{Z}^{n-1}}[e^{\tau\psi_{\varepsilon,g}(\cdot,0)}\theta_1\eta_{g,\mu}]^2_{1/2,\mathbb{R}^{n-1}}\\
\leq&C[e^{\tau\psi_{\varepsilon}(\cdot,0)}\theta_1]^2_{1/2,\mathbb{R}^{n-1}}
+C\mu\int_{\mathbb{R}^{n-1}}|\theta_1|^2e^{2\tau\psi_{\varepsilon}(x,0)}dx.
\end{aligned}
\end{equation}
Again by \eqref{7.191}, \eqref{6.26} and \eqref{6.20} we get
\begin{equation}\label{8.16}
\begin{aligned}
&\sum_{g\in \mathbb{Z}^{n-1}}[e^{\tau\psi_{\varepsilon,g}(\cdot,0)}J^{(1)}_{g,\mu}]^2_{1/2,\mathbb{R}^{n-1}}\\
\leq&C\left(\mu^3\sum_{\pm}\int_{\mathbb{R}^{n-1}}|w_{\pm}(x,0)|^2e^{2\tau\psi_{\varepsilon}(x,0)}dx
+\mu^2\sum_{\pm}[e^{\tau\psi_{\varepsilon}(\cdot,0)}w_{\pm}(\cdot,0)]^2_{1/2,\mathbb{R}^{n-1}}\right).
\end{aligned}
\end{equation}

We now go to the next term $\sum_{g\in \mathbb{Z}^{n-1}}[e^{\tau\psi_{\varepsilon,g}(\cdot,0)}J^{(2)}_{g,\mu}]^2_{1/2,\mathbb{R}^{n-1}}$. By \eqref{7.191}, \eqref{6.26} and \eqref{7.23}, we have that
\begin{equation}\label{8.19}
\begin{aligned}
&\sum_{g\in \mathbb{Z}^{n-1}}[e^{\tau\psi_{\varepsilon,g}(\cdot,0)}J^{(2)}_{g,\mu}]^2_{1/2,\mathbb{R}^{n-1}}\\
\leq&C\sum_{\pm}\left(\sum_{g\in \mathbb{Z}^{n-1}}[e^{\tau\psi_{\varepsilon}(\cdot,0)}\eta_{g,\mu}(A_{\pm}(\delta x_g,0)-A_{\pm}(\delta x,0))\nabla_{x,y} w_{\pm}\cdot \nu ]^2_{1/2,Q_{2/\mu}(x_g)}\right.\\
&\left.+\mu\sum_{g\in \mathbb{Z}^{n-1}}\int_{Q_{2/\mu}(x_g)}|A_{\pm}(\delta x_g,0)-A_{\pm}(\delta x,0)|^2|\nabla_{x,y} w_{\pm}|^2e^{2\tau\psi_{\varepsilon}(x,0)}dx\right)\\
\leq&C\sum_{\pm}\left(\delta^{2}\mu^{-2}\sum_{g\in \mathbb{Z}^{n-1}}[e^{\tau\psi_{\varepsilon}(\cdot,0)}\eta_{g,\mu}\nabla_{x,y} w_{\pm}(\cdot,0)]^2_{1/2,Q_{2/\mu}(x_g)}\right.\\
&\left.+\delta^{2}\mu^{-1}\int_{\mathbb{R}^{n-1}}|\nabla_{x,y}(w_{\pm}e^{\tau\psi_{\varepsilon}})(x,0)|^2dx+\delta^{2}\mu^{-1}\tau^2\int_{\mathbb{R}^{n-1}}|w_{\pm}(x,0)|^2e^{2\tau\psi_{\varepsilon}(x,0)}dx\right)
\\
\leq&C\sum_{\pm}\left(\delta^{2}\mu^{-2}[\nabla_{x,y}(w_{\pm}e^{\tau\psi_{\varepsilon}(\cdot,0)})]^2_{1/2,\mathbb{R}^{n-1}}+\delta^{2}\mu^{-2}\tau^2[e^{\tau\psi_{\varepsilon}(\cdot,0)}w_{\pm}]^2_{1/2,\mathbb{R}^{n-1}}\right.\\
&\left.+\delta^{2}\mu^{-1}\int_{\mathbb{R}^{n-1}}|\nabla_{x,y}w_{\pm}(x,0)|^2e^{2\tau\psi_{\varepsilon}(x,0)}dx+\delta^{2}\mu^{-1}\tau^2\int_{\mathbb{R}^{n-1}}|w_{\pm}(x,0)|^2e^{2\tau\psi_{\varepsilon}(x,0)}dx\right).
\end{aligned}
\end{equation}

Now we estimate $\sum_{g\in \mathbb{Z}^{n-1}}[e^{\tau\psi_{\varepsilon,g}(\cdot,0)}J^{(3)}_{g,\mu}]^2_{1/2,\mathbb{R}^{n-1}}$. In view of \eqref{7.191}, \eqref{6.26}, and \eqref{6.20}, we can obtain that
\begin{equation}\label{8.20}
\begin{aligned}
&\sum_{g\in \mathbb{Z}^{n-1}}[e^{\tau\psi_{\varepsilon,g}(\cdot,0)}J^{(3)}_{g,\mu}]^2_{1/2,\mathbb{R}^{n-1}}\\
\leq&C\left(\delta^{2}\sum_{\pm}[e^{\tau\psi_{\varepsilon}(\cdot,0)}w_{\pm}(\cdot,0)]^2_{1/2,\mathbb{R}^{n-1}}+\delta^{2}\mu\sum_{\pm}\int_{\mathbb{R}^{n-1}}|w_{\pm}(x,0)|^2e^{2\tau\psi_{\varepsilon}(x,0)}dx\right),
\end{aligned}
\end{equation}
where $C$ depends on $\lambda_0$, $M_0$ and $n$.
Finally, combining \eqref{8.12}, \eqref{8.13}, \eqref{8.16}, \eqref{8.19}, and \eqref{8.20} implies
\begin{equation}\label{8.21}
\sum_{g\in \mathbb{Z}^{n-1}}d^{(3)}_{g,\mu}\leq C\left(\tau\int_{\mathbb{R}^{n-1}}|\theta_1|^2e^{2\tau\psi_{\varepsilon}(x,0)}dx+[e^{\tau\psi_{\varepsilon}(\cdot,0)}\theta_1]^2_{1/2,\mathbb{R}^{n-1}}+R_4\right),
\end{equation}
where
\begin{equation*}
\begin{aligned}
R_4=&\delta^{2}\mu^{-2}\sum_{\pm}[\nabla_{x,y}(w_{\pm}e^{\tau\psi_{\varepsilon}})(\cdot,0)]^2_{1/2,\mathbb{R}^{n-1}}
+(\mu^2+\delta^{2}\mu^{-2}\tau^2)\sum_{\pm}[e^{\tau\psi_{\varepsilon}(\cdot,0)}w_{\pm}(\cdot,0)]^2_{1/2,\mathbb{R}^{n-1}}\\
&+\delta^2\varepsilon^{-1}\sum_{\pm}\int_{\mathbb{R}^{n-1}}|\nabla_{x,y}w_{\pm}(x,0)|^2e^{2\tau\psi_{\varepsilon}(x,0)}dx\\
&+(\varepsilon\tau^{2}+\delta^{2}\tau+\delta^{2}\mu^{-1}\tau^2)\sum_{\pm}\int_{\mathbb{R}^{n-1}}|w_{\pm}(x,0)|^2e^{2\tau\psi_{\varepsilon}(x,0)}dx.
\end{aligned}
\end{equation*}
Consequently, we have from \eqref{7.29}, \eqref{8.1}, \eqref{8.3}, \eqref{8.9} and \eqref{8.21} that
\begin{equation}\label{8.22}
\begin{aligned}
\Xi(w)\leq &C\left(\sum_{\pm}\int_{\mathbb{R}^n_{\pm}}|L_\delta(x,y,\partial)(w_\pm)|^2\,e^{2\tau\psi_{\varepsilon,\pm}}dxdy+[e^{\tau\psi_{\varepsilon}(\cdot,0)}\theta_1]^2_{1/2,\mathbb{R}^{n-1}}\right.\\
&+[\nabla_x(e^{\tau\psi_{\varepsilon}}\theta_0)(\cdot,0)]^2_{1/2,\mathbb{R}^{n-1}}+\tau^{3}\int_{\mathbb{R}^{n-1}}e^{2\tau\psi_{\varepsilon}(x,0)}|\theta_0(x)|^2dx\\
&\left.+\tau\int_{\mathbb{R}^{n-1}}e^{2\tau\psi_{\varepsilon}(x,0)}|\theta_1(x)|^2dx+R_5\right),
\end{aligned}
\end{equation}
where
\begin{equation*}
\begin{aligned}
R_5=
& \delta^2\mu^{-2}\sum_{\pm}\int_{\mathbb{R}^n_{\pm}}|D^2w_{\pm}|^2\,e^{2\tau\psi_{\varepsilon}}dxdy+\mu^{2}\sum_{\pm}\int_{\mathbb{R}^n_{\pm}}|Dw_{\pm}|^2\,e^{2\tau\psi_{\varepsilon}}dxdy\\
&+\mu^{4}\sum_{\pm}\int_{\mathbb{R}^n_{\pm}}|w_{\pm}|^2\,e^{2\tau\psi_{\varepsilon,\pm}}dxdy+(\mu+\delta^2\eps^{-1})\sum_{\pm}\int_{\mathbb{R}^{n-1}}|D{w}_{\pm}(x,0)|^2e^{2\tau\psi_{\varepsilon}(x,0)}dx\\
&+\mu\tau^2\sum_{\pm}\int_{\mathbb{R}^{n-1}}|w_{\pm}(x,0)|^2e^{2\tau\psi_{\varepsilon}(x,0)}dx+(\mu^4+\delta^{2}\mu^{-2}\tau^2)\sum_{\pm}[e^{\tau\psi_{\varepsilon}(\cdot,0)}w_{\pm}(\cdot,0)]^2_{1/2,\mathbb{R}^{n-1}}\\
&+\delta^{2}\mu^{-2}\sum_{\pm}[D(w_{\pm}e^{\tau\psi_{\varepsilon}})(\cdot,0)]^2_{1/2,\mathbb{R}^{n-1}}.
\end{aligned}
\end{equation*}

We now set $\delta=\varepsilon$ and choose a sufficiently small $\delta_0$  and a sufficiently large $\tau_0$, both depending on $\lambda_0,M_0,n$  such that if $\varepsilon\leq\delta_0$ and $\tau\geq\tau_0$, then $R_5$ on the right hand side of \eqref{8.22} can be absorbed by $\Xi(w)$ (defined in \eqref{l.s.}). In other words, we have proved that
\begin{equation}\label{8.23}
\begin{aligned}
&\sum_{\pm}\sum_{k=0}^2\tau^{3-2k}\int_{\mathbb{R}^n_{\pm}}|D^k{w}_{\pm}|^2e^{2\tau\psi_{\varepsilon}}dxdy+\sum_{\pm}\sum_{k=0}^1\tau^{3-2k}\int_{\mathbb{R}^{n-1}}|D^k{w}_{\pm}(x,0)|^2e^{2\psi_{\varepsilon}(x,0)}dx\\
&+\sum_{\pm}\tau^2[e^{\tau\psi_{\varepsilon}}w_{\pm}(\cdot,0)]^2_{1/2,\mathbb{R}^{n-1}}+\sum_{\pm}[\partial_y(e^{\tau\psi_{\varepsilon,\pm}}w_{\pm})(\cdot,0)]^2_{1/2,\mathbb{R}^{n-1}}
+\sum_{\pm}[\nabla_x(e^{\tau\psi_{\varepsilon}}w_{\pm})(\cdot,0)]^2_{1/2,\mathbb{R}^{n-1}}\\
\leq &C\left(\sum_{\pm}\int_{\mathbb{R}^n_{\pm}}|\mathcal{L}_\delta(x,y,\partial)(w_{\pm})|^2\,e^{2\tau\psi_{\varepsilon}}dxdy+[e^{\tau\psi_{\varepsilon}(\cdot,0)}\theta_1]^2_{1/2,\mathbb{R}^{n-1}}+[\nabla_x(e^{\tau\psi_{\varepsilon}}\theta_0)(\cdot,0)]^2_{1/2,\mathbb{R}^{n-1}}\right.\\
&\left.+\tau^{3}\int_{\mathbb{R}^{n-1}}|\theta_0|^2e^{2\tau\psi_{\varepsilon}(x,0)}dx+\tau\int_{\mathbb{R}^{n-1}}|\theta_1|^2e^{2\tau\psi_{\varepsilon}(x,0)}dx\right).
\end{aligned}
\end{equation}
Now, applying \eqref{8.23} to the function $w(x,y)=u(\delta x,\delta y)$ and by a standard change of variable, we have \eqref{8.24}, where $\phi_{\delta,\pm}$ is given by \eqref{wei}.

\end{document}